\documentclass[preprint,12pt]{elsarticle}

\usepackage{comment}
\usepackage{amssymb}
\usepackage{nomencl}
\usepackage{graphicx}%
\usepackage{multirow}%
\usepackage{amsmath,amssymb,amsfonts}%
\usepackage{amsthm}%
\usepackage{mathrsfs}%
\usepackage[title]{appendix}%
\usepackage{xcolor}%
\usepackage{textcomp}%
\usepackage{manyfoot}%
\usepackage{booktabs}%
\usepackage{algorithm}%
\usepackage{algorithmicx}%
\usepackage{algpseudocode}%
\usepackage{listings}%
\usepackage{tikz}
\usetikzlibrary{positioning}
\usepackage{xurl}
\usepackage{etoolbox}
\usepackage[margin = 3cm]{geometry}
\usepackage{accents}
\newcommand{\ubar}[1]{\underaccent{\bar}{#1}}

\newcommand{\blue}[1]{\textcolor{blue}{#1}}

\journal{Electric Power Systems Research}

\begin{document}

\begin{frontmatter}



\title{How Can Energy Communities Provide Grid Services? \\ A Dynamic Pricing Mechanism with Budget Balance, Individual Rationality, and Fair Allocation} 




\author[inst1]{Bennevis Crowley}

\affiliation[inst1]{organization={Department of Wind and Energy Systems, Technical University of Denmark},
            addressline={Elektrovej 325}, 
            city={Lyngby},
            postcode={2800}, 
            country={Denmark}}

\author[inst1]{Jalal Kazempour}
\author[inst1]{Lesia Mitridati}


\begin{abstract}
Following recent Danish legislation promoting energy communities, we explore how to enable these communities to provide grid services to distribution system operators. In particular, we focus on ``capacity limitation services", where we propose a bilateral agreement in which an energy community is given reduced grid import tariffs by setting a cap to its consumption level in certain hours. This requires a coordination mechanism between the community manager and the prosumers within the community. We enable this coordination by developing a bilevel optimization model to be solved by the community manager, aiming to set dynamic, i.e., time- and prosumer-differentiated, prices. This coordination mechanism enabled by dynamic pricing ensures desirable market properties including budget balance for the community manager and individual rationality for prosumers, while encouraging (but not guaranteeing) a fair allocation of collected benefits among prosumers.
\end{abstract}



\begin{keyword}
Energy communities \sep capacity limitation services \sep dynamic pricing \sep bilevel optimization \sep fair allocation
\end{keyword}

\end{frontmatter}


\newpage
\makenomenclature
\renewcommand\nomgroup[1]{%
  \item[\bfseries
  \ifstrequal{#1}{A}{Indexes and Sets}{%
  \ifstrequal{#1}{C}{Upper level variables}{%
  \ifstrequal{#1}{D}{Lower level variables}{%
  \ifstrequal{#1}{B}{Parameters}{%
  }}}}%
]}
\nomenclature[A]{$t \in \mathcal{T}$}{Set of time steps (hours)}
\nomenclature[A]{$n \in \mathcal{N}$}{Set of nodes of the grid inside the energy community}
\nomenclature[A]{$i \in \mathcal{I}$}{Set of energy community members (prosumers) }

\nomenclature[B]{$\lambda_t^{\mathrm{spot}}$}{Day-ahead (spot) price at time $t$ \ [DKK/kWh]}
\nomenclature[B]{$Y^\mathrm{im}$}{Distribution grid import tariff \ [DKK/kWh]}
\nomenclature[B]{$Y^\mathrm{ex}$}{Distribution grid export tariff \ [DKK/kWh]}
\nomenclature[B]{$\alpha^{\mathrm{DSO}}$}{Penalty for exceeding the capacity limitation \ [DKK/kWh]}
\nomenclature[B]{$\beta^{\mathrm{DSO}}$}{Discount on the grid import tariff \ [\%]}
\nomenclature[B]{$C^{\, \rm{ext}}_i$}{Baseline cost for community member $i$, calculated when it is not part of the community \ [DKK]}
\nomenclature[B]{$\alpha^{\mathrm{shed}}_{it}$}{Cost of load shedding for community member $i$ at time $t$ \ [DKK/kWh]}
\nomenclature[B]{$\bar{P}^{\mathrm{DSO}}_t$}{Capacity limitation imposed by the DSO at time $t$ \ [kW]}
\nomenclature[B]{$\bar{S}_n$}{Apparent power flow limit of the line upstream of node $n$ \ [p.u.]}
\nomenclature[B]{$R_n$}{Resistance of the line upstream of node $n$ \ [$\Omega$]}
\nomenclature[B]{$X_n$}{Reactance of the line upstream of node $n$ \ [$\Omega$]}
\nomenclature[B]{$\ubar{U}, \bar{U}$}{Lower and upper voltage limits at all nodes \ [p.u.] }
\nomenclature[B]{$PV_{it}$}{Photovoltaic power produced by community member $i$ at time $t$ \ [kW]}
\nomenclature[B]{$D_{it}$}{Demand of community member $i$ at time $t$ \ [kW]}
\nomenclature[B]{$\eta^{\rm{ch}}_i$}{Charging efficiency of the battery belonging to community member $i$ \ [\%]}
\nomenclature[B]{$\eta^{\rm{dis}}_i$}{Discharging efficiency of the battery belonging to community member $i$  \ [\%]}
\nomenclature[B]{$\bar{E}_i$}{Maximum state of charge of the battery belonging to community member $i$ \ [kWh]}
\nomenclature[B]{$\bar{P}^{\rm{bat}}_i$}{Maximum charging and discharging rate of the battery belonging to community member $i$ \ [kW]}
\nomenclature[B]{$\rho$}{Regularization constant for the maximum price in the objective function \ [-]}
\nomenclature[B]{$\gamma^{\mathrm{eq}}$}{Weight of the equal benefit distribution mechanism \ [-]}
\nomenclature[B]{$\gamma^{\mathrm{pro}}$}{Weight of the proportional benefit distribution mechanism \ [-]}
\nomenclature[B]{$S^{\mathrm{base}}$}{Base power of the distribution system \ [kVA]}

\nomenclature[C]{$p^{\mathrm{im}}_t$}{Imported active power from the upstream grid at time t \ [kW]}
\nomenclature[C]{$p^{\mathrm{ex}}_t$}{Exported active power to the upstream grid at time t \ [kW]}
\nomenclature[C]{$p^{\mathrm{pen}}_t$}{Amount of power exceeding the capacity limit at time t \ [kW]}
\nomenclature[C]{$x_{it}$}{Price set by the community manager for member $i$ at time $t$ \ [DKK/kWh]}
\nomenclature[C]{$w^+_i$}{Additional cost for community member $i$ compared to baseline \ [DKK]}
\nomenclature[C]{$w^-_i$}{Saved cost for community member $i$ compared to baseline \ [DKK]}
\nomenclature[C]{$v^+,v^-$}{SOS1 variables to ensure all community members benefit or lose \ [-]}
\nomenclature[C]{$\bar{x}$}{Maximum price across all hours and community members \ [DKK/kWh]}
\nomenclature[C]{$f^{\mathrm{p}}_{nt}$}{Active power flow entering node $n$ at time $t$ \ [kW]}
\nomenclature[C]{$f^{\mathrm{q}}_{nt}$}{Reactive power flow entering node $n$ at time $t$ \ [kVAR]}
\nomenclature[C]{$u_n$}{Squared voltage at node $n$ \ [p.u.]}
\nomenclature[C]{$q^{\mathrm{im}}_t$}{Imported reactive power from the upstream grid \ [kVAR]}
\nomenclature[C]{$q^{\mathrm{ex}}_t$}{Exported reactive power to the upstream grid \ [kVAR]}

\nomenclature[D]{$p^+_{it}$}{Imported active power by community member $i$ at time $t$ \ [kW]}
\nomenclature[D]{$p^-_{it}$}{Exported active power by community member $i$ at time $t$ \ [kW]}
\nomenclature[D]{$q^+_{it}$}{Imported reactive power by community member $i$ at time $t$ \ [kVAR]}
\nomenclature[D]{$q^-_{it}$}{Exported reactive power by community member $i$ at time $t$ \ [kVAR]}
\nomenclature[D]{$d^{\mathrm{shed}}_{it}$}{Amount of active power shed by community member $i$ at time $t$ \ [kW]}
\nomenclature[D]{$p^{\rm{ch}}_{it}$}{Battery charging by community member $i$ at time $t$ \ [kW]}
\nomenclature[D]{$p^{\rm{dis}}_{it}$}{Battery discharging by community member $i$ at time $t$ \ [kW]}
\nomenclature[D]{$e_{it}$}{State of charge of battery belonging to community member $i$ at time $t$ \ [kWh]}
\nomenclature[D]{$\lambda^{(.)}_{it}$}{Dual variable of equality constraints, where $(.)$ indexes the constraint. Time and consumer index vary based on the constraint \ [-]}
\nomenclature[D]{$\mu^{(.)}_{it}$}{Dual variable of inequality constraints, where $(.)$ indexes the constraint. Time and consumer index vary based on the constraint \ [-]}

\printnomenclature

\section{Introduction}\label{sec1}

The ongoing green transition is poised to revolutionize power systems, bringing new operational challenges and revolutionary ideas and opportunities. The growing concentration of distributed energy technologies and electrification of various sectors has led to increased risks related to line congestion, feeder overloading, thermal ratings of network components, voltage regulation, and power quality at all levels of the power system, especially the distribution grid. These challenges, among others, are outlined in more detail in \cite{CIGRE2014}. To mitigate these risks, distribution system operators (DSOs) must procure flexibility services from flexibility providers, such as aggregators and commercial end-users, among others, through various mechanisms. In the meantime, the growing active role of consumers and prosumers (consumers with local production) provides an opportunity to unlock flexibility within the distribution grid and address these challenges. As this transition progresses, new mechanisms and policies must come to the forefront to address these issues. In line with their aim to be a pioneering country in the green transition, the Danish Energy Agency released \textit{Market Model 3.0} \cite{Markedmodel2021} in 2021, a roadmap on how to make the Danish electricity market structure more flexible \cite{Gade2022}. A key focus of this roadmap is the development of new local flexibility products and markets that aim to meet congestion management needs throughout the distribution grid. Specifically, \cite{Markedmodel2021} highlights \textit{citizen energy communities} as a promising solution for gathering and coordinating flexibility from residential electricity users \cite{energycommunitiesMM32021}. Such an energy community is centered around flexible consumers and prosumers who are close to each other, mostly on the same feeder. As energy communities have the capability to access and control flexible assets from a variety of residential consumers and prosumers, they are well suited to offer flexibility services to DSOs.

In this context, this paper focuses on an energy community providing flexibility services to the DSO. The authors in \cite{Villar2018} provide a thorough overview of different solutions proposed in the literature for flexibility services in the distribution grid. The flexibility products and mechanisms mentioned in \cite{Villar2018} range from peer-to-peer flexibility markets to simple grid reinforcement. With many potential solutions being developed globally, this paper focuses only on flexibility services to limit congestion and feeder overload in the distribution grid. It will draw on established and recently proposed Danish policy and roadmaps to develop a novel local flexibility framework. Specifically, this paper focuses on the mechanisms the Danish Energy Agency identified in their technical report for Market Model 3.0 \cite{LocalFlexMM32021}.

\subsection{Literature review: Focus on capacity limitation services}
In \cite{LocalFlexMM32021}, flexibility mechanisms in distribution grids are split into two main categories: \textit{explicit} and \textit{implicit}. \textit{Implicit} flexibility mechanisms use a varying price signal to induce a change in demand or production, such as temporally and spatially differentiated tariffs. \textit{Explicit} flexibility mechanisms are those in which flexibility services, such as blackout contracts, special regulation services, or capacity limitation services, are traded between two parties, namely the DSO and a flexibility provider, through local flexibility markets or bilateral contracts. Amidst so many mechanisms, this paper focuses on \textit{capacity limitation services} traded through bilateral contracts between a DSO and an energy community, due to the following reasons.

Firstly, as the scope of this paper focuses on flexibility from residential end-users, blackout contracts are not a feasible long-term solution. This may be a solution for larger industrial-scale loads and is being implemented in Denmark \cite{StatusTarif2023}. However, this is not a desirable choice for an energy community that aims to provide security of supply to its members. Therefore, only solutions that circumvent delivery interruptions are considered. One possible mechanism is local flexibility markets, which are in use or under development in parts of Europe, such as the congestion management platform GOPACS in The Netherlands \cite{GOPACS}, the Piclo Flex Market in the UK \cite{Piclo2023} and NorFlex in southern Norway \cite{NODES2023}. These local markets allow for the submission of flexibility provision and purchasing bids and are subsequently cleared by a market operator, which may have to be the DSO if network data is required to clear the market \cite{Prat2023}. However, the added burden on DSOs to clear or even just submit bids in these markets remains a barrier to their implementation. Additionally, these markets rely on accurate consumer baseline consumption profiles, which come with their own complexities \cite{ziras2021}. Furthermore, special regulation services are already in use at the transmission level in Denmark \cite{EnerginetSystemydelser2023}, where the down-regulation of renewable generation is used to manage congestion within bidding zones. To the best of the authors' knowledge, an independently operated special regulation market for distribution grids does not exist, nor has it been the focus of research. Such a market would carry a similar, if not larger, burden as local flexibility markets and, as such, has not been considered further. Capacity limitation services, which are the focus of this paper, are one of the most feasible and practically straightforward flexibility services as they do not rely on complicated baseline consumption schedules for individual prosumers to deliver a congestion management service \cite{ziras2021}, and have as such been the subject of much research \cite{Ziras2017, Ziras2020, Heinrich2021, Chen2023}.

In a contract for capacity limitation services, an individual large consumer or a portfolio of flexible consumers or prosumers enters an agreement with the DSO to limit their power import from the grid. Such a contract can take different forms and vary in certain aspects, such as where in the grid and at what periods the limitation is enforced. Multiple papers have been summarized in Table \ref{tab:capacitylimitation} to give an overview of the current state-of-the-art of capacity limitation services.
\begin{table}[ht]
    \centering
    \begin{tabular}{|ccccc|}
         \hline 
         Reference          & Community        & Horizon  & Interaction with DSO   & Mechanism \\ \hline
         \cite{Ziras2017}   & No    & Daily         & Market-based  & Direct control             \\
         \cite{Ziras2020}   & No$^{*}$  & Monthly       & Market-based  & Direct control           \\
         \cite{Heinrich2021}& No   & Monthly       & Market-based  & Direct control            \\
         \cite{Chen2023}    & No$^{*}$ & Daily         & Market-based  & Direct control            \\ 
         This work          & Yes   & Daily         & Bilateral contract    & Indirect (price-based) \\                   \hline
    \end{tabular}
    \caption{Comparison of recent studies in capacity limitation services. $^{*}$Although not community-based, these studies consider aggregated electric vehicle loads.}
    \label{tab:capacitylimitation}
\end{table}
References \cite{Ziras2017} and \cite{Ziras2020} estimate the opportunity cost incurred by a portfolio of flexible assets and/or buildings as a consequence of providing capacity limitation services, which provides a basis to adequately price capacity limitation services before submitting the bids to a flexibility market. Reference \cite{Heinrich2021} then investigates how to design a market for a DSO to buy and price capacity limitation services. Lastly, \cite{Chen2023} models and characterizes the impact of these capacity limitation services on the network, offers valuable insights into system behavior, and empowers DSOs to request appropriate capacity limitation services better. All these studies consider capacity limitation services traded through a market platform. Organizing such a market platform may require an extensive administrative, communication, and computational burden from the DSO. In addition, the market-based proposals for local flexibility trading can be problematic due to the undesired strategic behavior of market participants, as discussed in \cite{utilities}. Since the applicability of market-based approaches seems too limited in practice, this paper focuses on bilateral capacity limitation contracts between the DSO and an energy community manager. Among others, a major challenge not yet addressed in the aforementioned studies is how individual assets can be coordinated and controlled by an aggregator or community manager to deliver on a signed capacity limitation contract. The community manager's main challenge is implementing an efficient strategy to control sufficient flexible assets to meet the DSO's demand. These studies assume direct control of the flexible assets. As direct control frameworks are limited by their scalability and the computation and information burden is placed on individual community members, \textit{indirect} control strategies, which carry less computational burden, may be preferred.

Historically, fixed-rate energy prices have been used for consumers throughout Europe. With increasing price volatility, many retailers have moved to time-of-use prices and other pricing mechanisms of a more dynamic nature to try and indirectly control demand. Furthermore, some stakeholders, such as the Danish DSOs, add additional layers of time-differentiated tariffs to encourage consumers to consider the timing of their energy usage \cite{Radius2023}. This naturally raises the question of whether or not consumers actually respond to price signals. In 2010, an analysis of household demand response programs conducted across the United States of America showed that simple dynamic price signals, such as time-of-use rates and critical-peak pricing, could elicit a demand response of up to 44\% when accompanied by enabling assets \cite{Faruqui2010}. Additionally, a recent increase of tariffs at peak demand hours has already shown to be effective, reducing peak loads in the Danish power system by as much as 10\% \cite{JanuarFebruar2023}. With more flexible assets and smart meters being integrated into the distribution grid, this number will only increase and the design of efficient price signals will become even more vital in harnessing demand-side flexibility. Current research efforts are investigating how to utilize the possibilities of dynamic pricing to elicit desirable behavior from residential end-users. For example, \cite{Grimm2021} studies the optimal design of retailer-prosumer tariffs while \cite{Askeland2023} investigates the optimal design of grid tariffs by the DSO to encourage flexible energy consumption.
Reference \cite{metric} proposes metrics to assess the performance of various grid tariff models. All these papers look at established pricing practices, such as price-of-use and capacity-based tariffs, instead of designing a novel pricing setup.

A potential approach to designing optimal tariffs is to allow for additional degrees of freedom, such as time and spatial differentiation. For example, \cite{Papavasiliou2018} introduces distribution locational marginal pricing, which differentiates prices in a distribution network by node through the use of Lagrangian multipliers, as is done in most wholesale markets. Another example is to implement stacked dynamic tariffs, where the grid tariff rates depend on the loading of the corresponding transformer, as discussed in \cite{EV}.
Another mechanism commonly used in pricing for energy systems is bilevel programming. Bilevel optimization models allow us to represent the underlying Stackelberg game structure which is present in many pricing problems. References \cite{Wogrin2020,pozo} outline the application of bilevel programs in the energy system, and include price signals as a potential application. Many papers go into greater detail and apply bilevel optimization models directly to demand response. References \cite{Meng2016,Alves2020,Venkatraman2022,Wei2022} propose different variations of bilevel optimization models for price-based demand response in which a retailer or utility determines price signals, seeking to maximize their profit while accounting for the cost-minimization problem of the following consumers.

\subsection{Contributions}

In view of the state-of-the-art and research gaps presented above, the contributions of this paper are threefold:

The first contribution is to design a novel coordination framework for delivering capacity limitation services from an energy community to the DSO, based on dynamic pricing. More specifically, a Stackelberg-game-inspired pricing mechanism is designed by which an energy community manager can coordinate the demand of the energy community members such that the resulting load from the community members stays below a limit set by the DSO. This approach aligns with recent developments in Danish energy policy, which allow for more differentiated prices and discounted grid tariffs for energy communities. This updated framework aims to facilitate the coordination between prosumers in an energy community and a DSO, enabling the provision of capacity limitation services. 
This coordination mechanism is formulated as a bilevel optimization problem, representing the interactions between the community manager and members, in which the community manager sets optimal prices which are both time- and space-dynamic while ensuring desirable economic properties including individual rationality and budget balance. The proposed setup uses an indirect control strategy within the energy community to provide explicit flexibility services to the DSO. This framework is then tested in a numerical study, in which we aim to provide a benchmark for how well an energy community can meet the requirements of a capacity limitation contract.

This paper's second contribution is an empirical analysis of the impact of the parameters stipulated in a capacity limitation contract on the community. An extensive numerical analysis of several possible combinations of capacity limitation variation and tariff discount parameters is conducted to identify any possible downsides and benefits for the energy community members. This analysis provides valuable insights for both energy communities and DSOs regarding what kind of parameters could be mutually beneficial when designing such a contract in the future.

Lastly, the third contribution of this paper is to propose a straightforward mechanism for the distribution of community benefits among its members. For this purpose, two regularization terms that can be added to the objective function of the community manager are introduced. Each regularization term corresponds to a different interpretation of ``fair" allocation, namely \textit{equal} and \textit{proportional} distribution. This regularization technique enables the community manager to systemically consider the distribution of benefits among the members when setting optimal prices. These regularization terms are then compared with respect to their impact on the benefits earned by each community member.

The remainder of this paper is structured in the following manner: Section \ref{sec:framework} introduces the proposed framework and coordination method to deliver capacity limitation services to the DSO, Section \ref{sec:models} describes the developed bilevel model, Section \ref{sec:results} shows the results to a case study, and Section \ref{sec:conclusion} summarizes the conclusions drawn and outlines important future works.

\section{Proposed Framework} \label{sec:framework}

This section explains the design of the envisioned capacity limitation services and provides insight into the daily workings of the proposed coordination framework for all involved parties, including the DSO, community manager, and community members.

\subsection{Capacity limitation service between the DSO and energy community} \label{sec:coordination}

To clarify the envisioned capacity limitation services, we provide an illustrative example of an energy community and its connection to the larger grid.  Figure \ref{fig:example} shows 6 residential households\footnote{This number is arbitrary and kept small to keep the illustration simple.} in a purple box (left-hand side) that is connected to the larger grid through a single connection point. 
The community delivers the capacity limitation service on this single connection point, indicated by the thick red bars around the line. The gray dashed box in Figure \ref{fig:example} contains a plot of the capacity limitation service in action. The typical consumption of the community is shown in the dashed black line, the capacity limit in the red horizontal line, and the resulting shifted demand profile in blue. The adjusted profile adheres to the capacity limit by shifting some consumption to earlier and later in the day when the capacity limit is not a limiting factor for the community import. Note that the capacity limitation is not necessarily a physical limitation of the system but rather an agreed-upon value between the DSO and the energy community.

\pgfdeclarelayer{nodelayer}
\pgfdeclarelayer{edgelayer}
\pgfdeclarelayer{background}

\tikzstyle{none}=[inner sep=0pt]
\tikzstyle{Time axis}=[|->, line width = 1pt]
\tikzstyle{DSO edge}=[-, draw=black, fill={rgb,255: red,153; green,255; blue,255}]
\tikzstyle{arrow}=[->,line width = 0.6pt]
\tikzstyle{scope}=[-, dashed, fill={rgb,255: red,191; green,191; blue,191}]
\tikzstyle{ECM edge}=[-, dashed, fill={rgb,255: red,177; green,156; blue,217}]
\tikzstyle{Consumer edge}=[-, fill={rgb,255: red,153; green,255; blue,153}]
\tikzstyle{edge}=[-, line width = 1pt]
\tikzstyle{grid}=[circle, fill = none, draw = black, minimum size = 0.5cm, line width = 1pt]
\tikzstyle{reference node}=[circle, fill = white, draw = black, minimum size = 0.75cm]
\tikzstyle{node}=[circle, fill = {rgb,255: red,153; green,255; blue,153}, draw = black, minimum size = 0.75cm]
\tikzstyle{blue arrow}=[->, draw = blue, line width = 2pt]
\tikzstyle{new edge style 0}=[-, draw = blue, line width = 1pt]
\tikzstyle{limit}=[-, draw = red, line width = 2pt]
\tikzstyle{point}=[coordinate]

\begin{figure} [t]
    \centering
    \begin{tikzpicture}
    \pgfsetlayers{background,nodelayer}
    \begin{pgfonlayer}{nodelayer}
        \node[none] (grid) at (0,0) {\includegraphics[width=1cm]{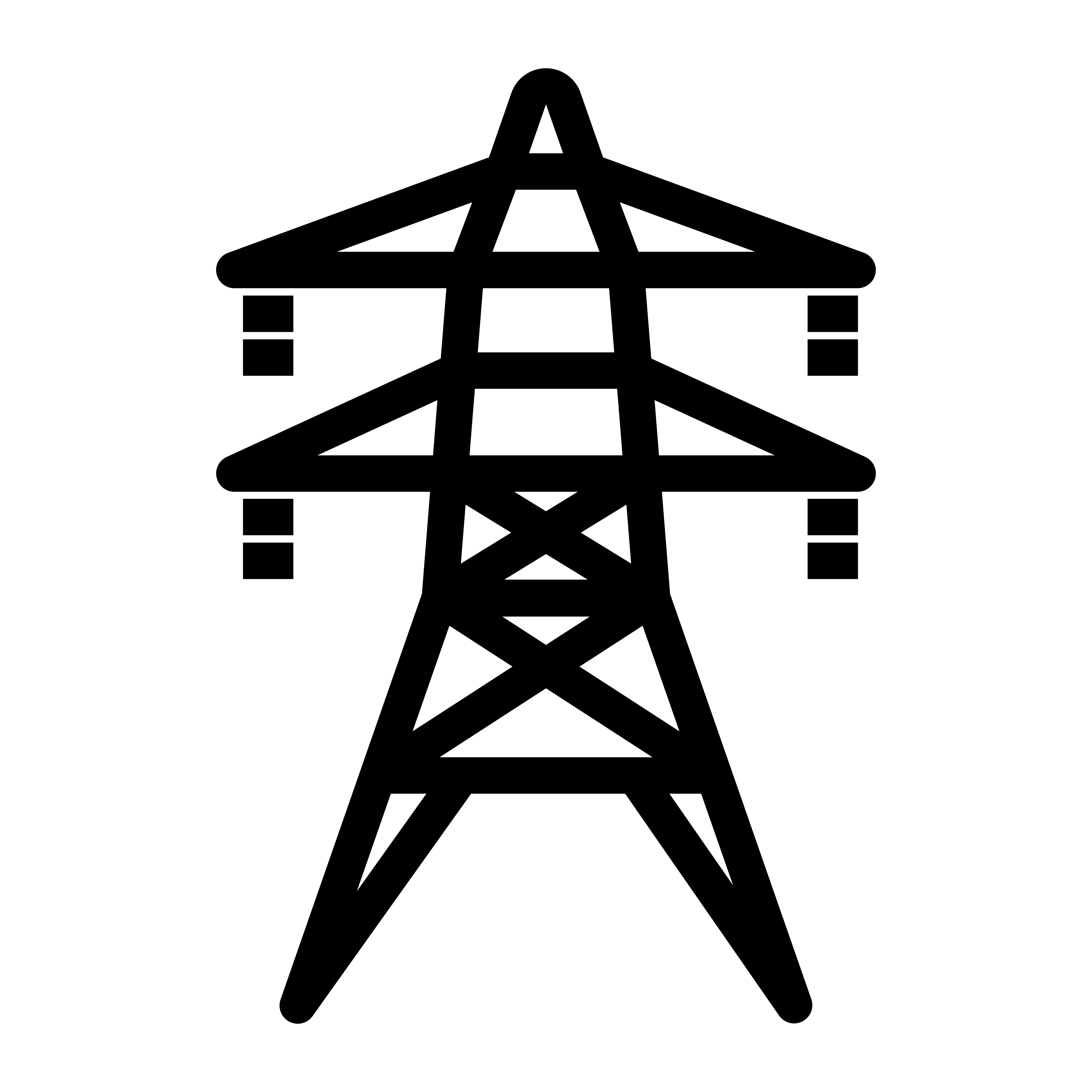}};
        \node[point, below = 1 of grid] (meeting point) {};
        \draw[edge] (grid) -- (meeting point);
        \node[point, left = 4 of meeting point] (EC1) {};
        \draw[edge] (meeting point) -- (EC1);
        \node[point, left = 1 of EC1] (EC2) {};
        \draw[edge] (EC1) -- (EC2);
        \node[point, left = 1 of EC2] (EC3) {};
        \draw[edge] (EC2) -- (EC3);
        \node[none, above = 0.5 of EC1] (House1) {{\includegraphics[width=0.75cm]{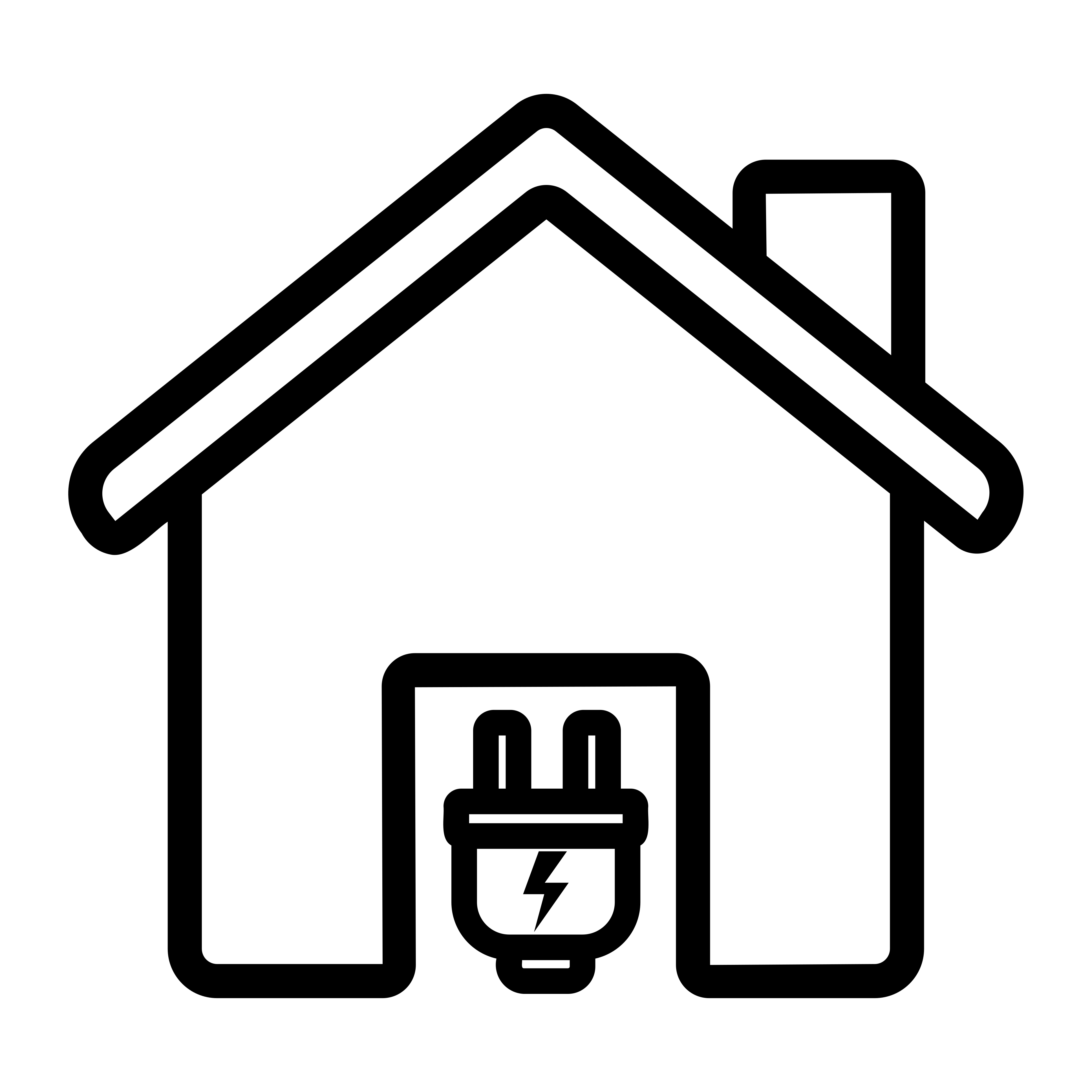}}};
        \draw[edge] (EC1) -- (House1);
        \node[none, above = 0.5 of EC2] (House2) {{\includegraphics[width=0.75cm]{haus.png}}};
        \draw[edge] (EC2) -- (House2);
        \node[none, above = 0.5 of EC3] (House3) {{\includegraphics[width=0.75cm]{haus.png}}};
        \draw[edge] (EC3) -- (House3);
        \node[none, below = 0.5 of EC3] (House4) {{\includegraphics[width=0.75cm]{haus.png}}};
        \draw[edge] (EC3) -- (House4);
        \node[none, below = 0.5 of EC2] (House5) {{\includegraphics[width=0.75cm]{haus.png}}};
        \draw[edge] (EC2) -- (House5);
        \node[none, below = 0.5 of EC1] (House6) {{\includegraphics[width=0.75cm]{haus.png}}};
        \draw[edge] (EC1) -- (House6);
        \node[none, above = 0.2 of House2](EClabel) {Energy Community};
        \node[point,above left = 0.1 and 1.5 of meeting point] (limit1){};
        \node[point,left = 0.75 of limit1] (limit2){};
        \draw[limit] (limit1) -- (limit2);
        \node[point, below left = 0.1 and 1.5 of meeting point] (limit3){};
        \node[point, left = 0.75 of limit3] (limit4){};
        \draw[limit] (limit3) -- (limit4);
        \node[circle, left = 1.35 of meeting point, minimum size = 1cm, draw = gray,line width=1,dashed](zoom start){};
        \node[below = 1 of meeting point, draw = gray, dashed, line width = 1] (zoom graph) {\includegraphics[width = 4cm]{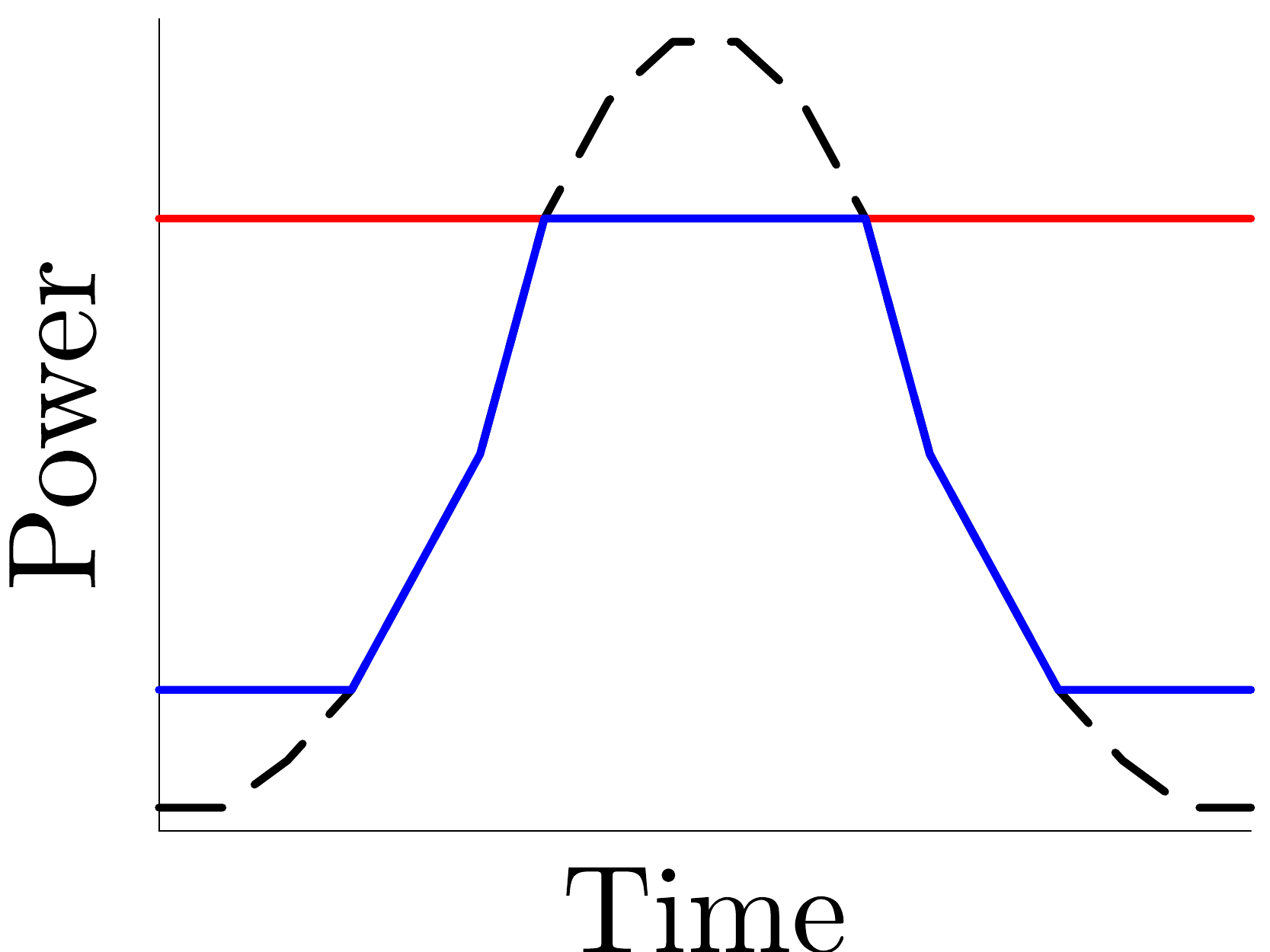}};
        \draw[dashed, draw=gray,line width = 1] (zoom start) -- (zoom graph);
        \node[point, right = 4 of meeting point] (grid1) {};
        \draw[edge] (meeting point) -- (grid1){};
        \node[none, above = 0.5 of grid1] (ext house) {\includegraphics[width=0.75cm]{haus.png}};
        \draw[edge] (grid1) -- (ext house){};
        \node[none, right = 0.5 of grid1] (right end) {\hspace{2mm}$\hdots$};
        \draw[edge] (grid1) -- (right end){};
        \node[none, below = 0.5 of grid1] (bottom end){$\vdots$};
        \draw[edge] (grid1) -- (bottom end){};
    \end{pgfonlayer}
    \begin{pgfonlayer}{background}
        \node[none, above right = 0.1 of House1] (community1) {};
        \node[none, above left = 0.1 of House3] (community2) {};
        \node[none, below left = 0.1 of House4] (community3) {};
        \node[none, below right = 0.1 of House6] (community4) {};
        \draw [style=ECM edge] (community1.center) to (community2.center) to (community3.center)to (community4.center) to cycle;
    \end{pgfonlayer}
    \end{tikzpicture}
    \caption{Illustration of capacity limitation (red lines) being imposed on an energy community (purple box).}
    \label{fig:example}
\end{figure}

To incentivize the energy community to adhere to a capacity limit, the DSO must reward the community for agreeing to follow the limit. The proposed framework proposes that the DSO leverages grid tariffs to incentivize residential prosumers to provide flexibility services. This framework is based on the idea that in exchange for providing the capacity limitation service at their connection point to the DSO-operated grid, the energy community receives a \textit{discount} on the grid tariffs for all power flows that occur inside the community, i.e., downstream from the line where the capacity limit is imposed. Leveraging tariffs as an incentive naturally depends on tariffs being a sufficiently large portion of electricity costs. This is currently the case in Denmark, as current Danish transport tariffs make up more than 20\% of household electricity bills over a year \cite{2024Elprisen}, and in winter, this percentage can regularly reach up to 50\% in peak hours \cite{2023TarifferNetabonnement}. In addition to providing the benefit of a discount on internal power flows, this framework also introduces a kW-based penalty imposed by the DSO that is charged to the community for any consumption exceeding the agreed-upon capacity limit.

\subsection{Coordination between community manager and members}

The community manager is responsible for coordinating the energy consumption/production of the community members to meet the agreed-upon capacity limitation with the DSO and share the benefit among all prosumers. The proposed framework achieves this by sending dynamic price signals $x_{it}$ which cover electricity costs and grid tariffs. These prices are temporally and spatially dynamic, meaning that each community member $i \in \mathcal{I}$ receives a different price for all time steps $t \in \mathcal{T}$ of the next day. Each community member then independently optimizes their own electricity consumption/production $p^+_{it}/p^-_{it}$ based on the price signals $x_{it}$ to minimize their own energy costs. Community members respond directly to individualized prices only specified within the community, not the day-ahead prices and tariffs of the larger grid.

\begin{figure*}
    \pgfsetlayers{background,edgelayer,nodelayer}
    \centering
    \begin{tikzpicture}[scale = 0.35]
    	\begin{pgfonlayer}{nodelayer}
    		\node [style=none] (0) at (-18, 26) {};
    		\node [style=none] (1) at (-14, 26) {};
    		\node [style=none] (2) at (-14, 22) {};
    		\node [style=none] (4) at (-18, 22) {};
    		\node [style=none] (6) at (-16, 24) {DSO};
    		\node [style=none] (7) at (-4, 26) {};
    		\node [style=none] (9) at (4, 22) {};
    		\node [style=none] (10) at (4, 26) {};
    		\node [style=none] (11) at (-4, 22) {};
    		\node [style=none] (15) at (0, 24.75) {Community};
    		\node [style=none] (16) at (0, 23.25) {Manager};
    		\node [style=none] (17) at (-14, 24) {};
    		\node [style=none] (19) at (-7.5, 24) {};
    		\node [style=none] (20) at (12, 26) {};
    		\node [style=none] (21) at (12, 22) {};
    		\node [style=none] (22) at (18, 22) {};
    		\node [style=none] (23) at (18, 26) {};
    		\node [style=none] (24) at (15, 24) {Prosumers};
    		\node [style=none] (25) at (7.5, 24) {};
    		\node [style=none] (26) at (12, 24) {};
    		\node [style=none] (27) at (-20, 10.5) {};
    		\node [style=none] (28) at (18, 10.5) {};
    		\node [style=none] (33) at (-20, 16) {};
    		\node [style=none] (34) at (-20, 12) {14:00};
    		\node [style=none] (39) at (2, 16) {};
    		\node [style=none] (40) at (2, 14) {};
    		\node [style=none] (41) at (7, 14) {};
    		\node [style=none] (42) at (7, 16) {};
    		\node [style=none] (43) at (-2, 16) {};
    		\node [style=none] (44) at (-2, 14) {};
    		\node [style=none] (45) at (-7, 14) {};
    		\node [style=none] (46) at (-7, 16) {};
    		\node [style=none] (47) at (-0.75, 22) {};
    		\node [style=none] (52) at (-4, 17) {};
    		\node [style=none] (53) at (0, 17) {};
    		\node [style=none] (54) at (4, 17) {};
    		\node [style=none] (55) at (-4, 17) {};
    		\node [style=none] (56) at (-4.5, 15) {\tiny{Prosumer 1}};
    		\node [style=none] (57) at (0, 15) {\huge{...}};
    		\node [style=none] (58) at (4.5, 15) {\tiny{Prosumer n}};
    		\node [style=none] (59) at (-4.75, 16) {};
    		\node [style=none] (61) at (4.75, 16) {};
    		\node [style=none] (62) at (0.75, 22) {};
    		\node [style=none] (63) at (1.25, 22) {};
    		\node [style=none] (64) at (-1.25, 22) {};
    		\node [style=none] (65) at (-4.25, 16) {};
    		\node [style=none] (66) at (4.25, 16) {};
    		\node [style=none] (71) at (-5.5, 17) {\small{$\mathbf{x}$}};
    		\node [style=none] (73) at (2.5, 17) {\small{$\mathbf{x}$}};
                \node [style=none] (74) at (-2.25, 17) {\small{$\mathbf{p,q}$}};
                \node [style=none] (77) at (5.75, 17) {\small{$\mathbf{p,q}$}};
    		\node [style=none] (78) at (-16, 11) {};
    		\node [style=none] (79) at (-16, 10) {};
    		\node [style=none] (80) at (0, 11) {};
    		\node [style=none] (81) at (0, 10) {};
    		\node [style=none] (82) at (-16, 9) {Step 2};
                \node [style=none] (83) at (-10.75, 25) {\small{$\beta_t^{\rm{DSO}},\alpha_t^{\rm{DSO}}$}};
                \node [style=none] (90) at (-10.75,23) {\small{$\bar{P}^{\rm{DSO}}_t$}};
    		\node [style=none] (84) at (-10.75, 21.25) {Step 3};
    		\node [style=none] (85) at (0, 9) {Step 4};
    		\node [style=none] (86) at (15, 11) {};
    		\node [style=none] (87) at (15, 10) {};
    		\node [style=none] (89) at (15, 9) {Step 6};
    		\node [style=none] (94) at (-7, 23.5) {};
                \node [style=none] (99) at (9.75, 25) {\small{$x_{it}$}};
    		\node [style=none] (100) at (9.75, 23) {Step 5};
    		\node [style=none] (101) at (7, 23.5) {};
    		\node [style=none] (104) at (-20, 9) {Step 1};
    		\node [style=none] (108) at (7.5, 26.5) {};
    		\node [style=none] (109) at (7.5, 13.5) {};
    		\node [style=none] (110) at (-7.5, 13.5) {};
    		\node [style=none] (111) at (-7.5, 26.5) {};
                \node [style=none] (112) at (-16, 12) {$\approx$14:15};
                \node [style=none] (113) at (0, 12) {$\approx$14:30};
                \node [style=none] (114) at (15, 12) {$\approx$14:30+};
    	\end{pgfonlayer}
    	\begin{pgfonlayer}{edgelayer}
    		\draw [style=DSO edge] (2.center)
    			 to (1.center)
    			 to (0.center)
    			 to (4.center)
    			 to cycle;
    		\draw [style=ECM edge] (11.center)
    			 to (9.center)
    			 to (10.center)
    			 to (7.center)
    			 to cycle;
    		\draw [style=arrow] (17.center) to (19.center);
    		\draw [style=Consumer edge] (21.center)
    			 to (22.center)
    			 to (23.center)
    			 to (20.center)
    			 to cycle;
    		\draw [style=arrow] (25.center) to (26.center);
    		\draw [style=Time axis] (27.center) to (28.center);
    		\draw [style=Consumer edge] (43.center)
    			 to (46.center)
    			 to (45.center)
    			 to (44.center)
    			 to cycle;
    		\draw [style=Consumer edge] (39.center)
    			 to (40.center)
    			 to (41.center)
    			 to (42.center)
    			 to cycle;
    		\draw [style=arrow] (64.center) to (59.center);
    		\draw [style=arrow] (65.center) to (47.center);
    		\draw [style=arrow] (62.center) to (66.center);
    		\draw [style=arrow] (61.center) to (63.center);
    		\draw [style=edge] (78.center) to (79.center);
    		\draw [style=edge] (80.center) to (81.center);
    		\draw [style=edge] (86.center) to (87.center);
    	\end{pgfonlayer}
            \begin{pgfonlayer}{background}
                \draw [style=scope] (109.center)
    			 to (110.center)
    			 to (111.center)
    			 to (108.center)
    			 to cycle;
            \end{pgfonlayer}
    \end{tikzpicture}
    \caption{The decision-making and information exchange sequence among different stakeholders in the proposed framework.}
    \label{fig:illustration}
\end{figure*}
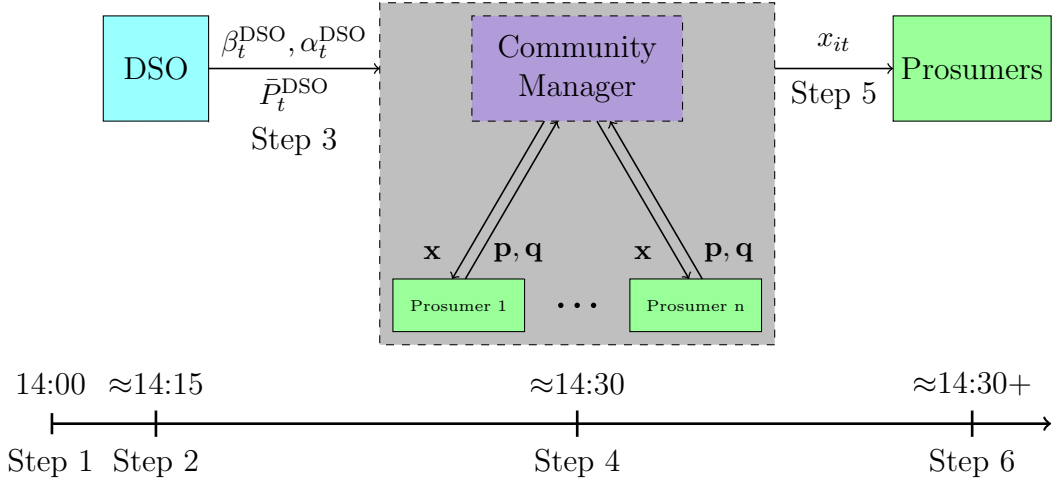

The timeline of decision-making and information exchange among stakeholders is illustrated in Figure \ref{fig:illustration} and summarized in the following steps:
\begin{enumerate}
    \item The day-ahead electricity market prices $\lambda^{\rm{spot}}_t$ for each hour of the next day are released by the wholesale market operator (e.g., at 14:00 by Nord Pool).
    \item The DSO identifies the expected needs and costs of congestion management services in the distribution grid and calculates the required capacity limitation.
    \item The DSO communicates the required capacity limitations $\bar{P}_t^{\rm{DSO}}$ for each hour $t$ of the next day to the community manager, as well as penalties $\alpha_t^{\rm{DSO}}$ for violating these limitations and the offered grid tariff reduction $\beta_t^{\rm{DSO}}$ for internal power flows within the community.
    \item Given the DSO's offer in Step 3, the community manager sets prices $x_{it}$ for each community member $i$ and hour $t$ of the next day.
    \item The community manager communicates the price signals $x_{it}$ to all community members $i$.
    \item Given the price signals $x_{it}$ each community member $i$ optimizes their own energy consumption/production profile $p^+_{it}/p^-_{it}$ for the next day.
\end{enumerate}

The main focus of this paper is on Step 4 (the gray box in Figure \ref{fig:illustration}), also known as the price-setting problem. Therefore, a mathematical program for the community manager's pricing problem is developed in Section \ref{sec:models}. The mechanism by which the DSO sets these capacity limitation contracts in Steps 2-3 (i.e., $\bar{P}_t^{\rm{DSO}}, \alpha_t^{\rm{DSO}}, \beta_t^{\rm{DSO}}$) is outside the scope of this paper. 

\subsection{Technical  requirements for practical implementation}

To ensure that the above-proposed framework is indeed implementable, certain technical capabilities would be required within the energy community. Firstly, communication and measurement infrastructure must be present within the community. This includes hardware requirements, such as a stable internet connection and sufficient metering infrastructure at the household level (a smart meter per community member). There are also certain software requirements that must be present, such as a secure two-way communication interface between the DSO and energy community manager (for the capacity limitation parameters) as well as between the energy community manager and the individual community members (for dissemination of the prices to consumers and timely collection of relevant parameters for the community manager).

Along with the technical requirements, some other aspects must be considered before implementing such a framework in practice. The optimization problem of the energy community manager must be solved in a reasonable time duration and not be computationally expensive or scale at an exponential rate, such that the prices can be shared with the community members in good time. Additionally, grid data must be shared with the energy community manager to ensure the pricing mechanism provides grid-aware prices. Furthermore, a registry of flexible assets in an energy community is required. Without this knowledge, it is not straightforward for the energy community manager to decide what flexibility services they can offer the DSO. All these requirements need data sharing, which inherently comes with privacy risks. Therefore, policies must be in place regarding how such data must be handled to ensure privacy for all involved parties. Without consideration of these aspects, which are outside this paper's scope, such a framework could prove difficult to implement in practice.

\section{Proposed Mathematical Model} \label{sec:models}

First, we introduce the game-theoretical background for the proposed dynamic pricing scheme. Second, we provide a detailed formulation of the proposed bilevel program representing the community manager's problem. Third, we discuss the distribution of community benefits among prosumers.  Finally, we explain the reformulation methods used to solve this mathematical optimization problem efficiently.

\subsection{The Stackelberg game} \label{sec:stackelberg}

A Stackelberg game classically consists of a single leading agent, who needs to make an initial decision that impacts the subsequent decisions of the following agent or agents, whose actions, in turn, impact the resulting outcome for the leader. This structure is represented in Steps 4-6 of the framework described in Section \ref{sec:framework}. The leader (community manager) optimizes their own decisions (price signals) first (Step 4), which then impacts (Step 5) the optimal decisions (consumption/production) of multiple followers (community members). In turn, the followers' reaction (Step 6) impacts the leader's optimal cost through the total power imported/exported from the upstream grid. Therefore, the community manager must anticipate the reaction of the community members to set the optimal price signals that will elicit the desired reaction.

To model this setup, we propose a formulation of the aforementioned community manager's problem as a bilevel program, in which a set of lower-level optimization problems, representing the optimal reaction of the community members, is embedded as constraints of the upper-level optimization problem, representing the optimal decisions of the community manager. A compact formulation of this bilevel program is given in \eqref{compact1}, including the upper-level problem \eqref{compact2}-\eqref{compact3}, and the set of lower-level problems \eqref{compact4}-\eqref{compact6}, one per community member $i$, as follows:
\begingroup
\allowdisplaybreaks
\begin{subequations} \label{compact1}
\begin{alignat}{1}
    &\underset{\underset{\textbf{p}^{\rm{im}},\textbf{p}^{\rm{ex}}}{\textbf{x}_i,}}{\operatorname{min}} \text{cost}^{\rm{C}}(\textbf{x}_i,\textbf{p}^{\rm{im}} ,\textbf{p}^{\rm{ex}},\textbf{p}_i) \label{compact2} \\ 
    & \hspace{0.25cm} \text{s.t.} \ g^{\rm{C}}(\textbf{x}_i,\textbf{p}^{\rm{im}} ,\textbf{p}^{\rm{ex}},\textbf{p}_i) = 0 \\
    & \hspace{0.25cm} \hspace{0.54cm} \ h^{\rm{C}}(\textbf{x}_i,\textbf{p}^{\rm{im}} ,\textbf{p}^{\rm{ex}},\textbf{p}_i) \leq 0 \label{compact3} \\ 
    & \hspace{0.25cm} \hspace{0.54cm} \ \textbf{p}_i \in \underset{\textbf{p}_i}{\text{argmin}} \ \Big\{ \text{cost}_i^{\rm{P}}(\textbf{x}_i,\textbf{p}_i) \label{compact4}\\
    & \hspace{0.25cm} \hspace{0.54cm} \hspace{2.12cm} \text{s.t.} \  g_i^{\rm{P}}(\textbf{p}_i) = 0 \label{compact5}\\
    & \hspace{0.25cm} \hspace{0.54cm} \hspace{2.12cm} \hspace{0.54cm} \ h_i^{\rm{P}}(\textbf{p}_i) \leq 0 \Big\} \ \forall i \in \mathcal{I}. \ \label{compact6}
\end{alignat}
\end{subequations}
\endgroup

This bilevel program is solved by the community manager, given the day-ahead electricity market prices $\lambda^{\rm{spot}}_t$, the grid tariffs $\textbf{Y}$, as well as the DSO capacity limitation contract parameters ($\bar{P}_t^{\rm{DSO}}, \beta_t^{\rm{DSO}}, \alpha_t^{\rm{DSO}}$). In the upper-level problem, the community manager aims at finding the optimal vectors of price signals $\textbf{x}_i$ which minimize the total community cost, i.e.,  $\text{cost}^{\rm{C}}(\textbf{x}_i,\textbf{p}^{\rm{im}},\textbf{p}^{\rm{ex}},\textbf{p}_i)$ in \eqref{compact2}. This objective function is subject to community-wide constraints related to limitations on the imported and exported power, 
feasibility of power flows within the community, and desirable economic properties such as budget balance and individual rationality. Additionally, the community manager anticipates how each community member $i$  optimizes their consumption/production profile $\textbf{p}_i$ in response to the price signals $\textbf{x}_i$. This optimal response of the community members is modeled as the solutions to a set of lower-level optimization problems, one per community member $i$, which are embedded as constraints in the upper-level problem. In each lower-level problem, a given community member $i$ minimizes their daily energy cost, i.e.,  ${\rm{cost}}_i^{\rm{P}}(.)$, subject to a set of internal constraints \eqref{compact5}-\eqref{compact6}, for the given price signals $\textbf{x}_i$ (treated as parameters in the lower-level problems). Implementing this dynamic pricing scheme requires the community manager to have perfect knowledge of the community members' parameters, which might raise privacy and communication challenges in practice. Additionally, in the case of a multiplicity of solutions to the lower-level problems, the upper-level problem will choose an ``optimistic" solution, which minimizes its own objective. However, in practice, the community manager cannot guarantee which optimal solutions individual members would choose. Therefore, we see the proposed framework as an \textit{ideal benchmark}, which provides an upper bound on how effective a dynamic pricing scheme can be to incentivize and coordinate prosumers in an energy community to procure capacity limitation services to the DSO.

\subsection{Detailed mathematical formulation} \label{subsec:detailed models}

\subsubsection{Upper-level problem (dynamic pricing)} 
The full upper-level problem is described in \eqref{eqs:upper level} and makes use of the variable set $\Omega = \{x_{it}$, $p^{\rm{im}}_t$, $p^{\rm{ex}}_t$, $p^{\rm{pen}}_t$, $q^{\rm{im}}_t$, $q^{\rm{ex}}_t$, $f^{\rm{q}}_{lt}$, $f^{\rm{p}}_{lt}$, $u_{nt}$, $\bar{x}$, $w^+_i$, $w^-_i$, $v^+$, $v^-$\}. In addition, $\mathcal{T}$ and $\mathcal{I}$ are sets for times and prosumers, respectively. The upper-level objective function \eqref{eq:upper_obj}, minimizing the total community cost, reads as
%
\begin{subequations} \label{eqs:upper level}
\begin{align}
    \nonumber \underset{\Omega}{\operatorname{min}} \quad \sum_{t \in \mathcal{T}} &\Bigg[ p^{\rm{im}}_t \Big(\lambda^{\rm{spot}}_t + Y^{\rm{im}}_t\Big) - p^{\rm{ex}}_t \Big(\lambda^{\rm{spot}}_t - Y^{\rm{ex}}_t\Big) \\
    \nonumber & + (1-\beta^{\rm{DSO}}_t) Y^{\rm{im}}_t \Big( \sum_{i \in \mathcal{I}} p^+_{it} - p^{\rm{im}}_t \Big) \\
    & + \alpha^{\rm{DSO}}_t p^{\rm{pen}}_t + \alpha^{\rm{shed}} d^{\, \rm{shed}}_{it} \Bigg] + \rho \bar{x}^2.  \label{eq:upper_obj}
\end{align}

The first line corresponds to energy trades and grid tariffs at the interface of the energy community and the DSO-operated upstream grid. The community pays (is paid) for the imported power $p^{\rm{im}}_t$ from (exported power $p^{\rm{ex}}_t$ to) the upstream grid at given market prices $\lambda^{\rm{spot}}_t$. Additionally, the community pays grid tariffs that are set by the DSO for both import ($Y^{\rm{im}}_t$) and export ($Y^{\rm{ex}}_t$) at the interface, i.e., for the usage of the feeder between the community and the upstream grid. 

The second line of \eqref{eq:upper_obj} pertains to discounted grid tariffs for power flows that stay within the radial grid of the community. The provided DSO discount is implemented by the coefficient  $\beta^{\rm{DSO}}_t$, whose value lies between zero and one. Therefore, $(1-\beta^{\rm{DSO}}_t)Y^{\rm{im}}_t$ is the discounted import tariff rate at time $t$ for the usage of lines within the community. Symbols $p^+_{it}$ and $p^-_{it}$ denote the power purchase and sale of prosumer $i$ at time  $t$, and are determined in the lower-level problem. 
Therefore, for time $t$, the term $\sum_{i}p^+_{it} - p^{\rm{im}}_t$ represents the difference between total purchased power and power imported from the grid, i.e., the internal power flow of the community.

The third line of \eqref{eq:upper_obj} has three cost elements. The first term accounts for the DSO penalty, such that the parameter $\alpha^{\rm{DSO}}_t$ is the penalty rate set by the DSO for every 1  kW power in the feeder at time $t$ that exceeds the cap $\bar{P}_t^{\rm{DSO}}$, where the variable $p^{\rm{pen}}_t$ is the amount of exceeded power. The second term penalizes any possible load shedding $d^{\, \rm{shed}}_{it}$ in the lower-level problem by prosumer $i$ at time $t$ at a predetermined value of lost load $\alpha^{\rm{shed}}$. This prevents the community manager from choosing a solution prioritizing load shedding over other solutions. Lastly, the third term is a quadratic regularization term, which aims to converge to a unique solution for the price setting. As there is no upper bound for the dynamic prices $x_{it}$, there could be multiple solutions for $x_{it}$, leading to the same value for the objective function of the community. This regularization term prioritizes solutions with the lowest maximum price across all prosumers, minimizing the squared maximum price observed $\bar{x}$, multiplied by a small positive weight $\rho$. 

The upper-level constraints, summarized in \eqref{compact3}, include community-level constraints \eqref{eq:budget_balance} to \eqref{eq:positive variables}. Constraint \eqref{eq:budget_balance} obtains the budget balance property for the energy community manager, meaning neither taking a loss nor making a profit:
\begin{align}
    \nonumber  \sum_{t \in \mathcal{T}} \sum_{i \in \mathcal{I}} \Big( p^+_{it} - p^-_{it} \Big) x_{it}  = &\sum_{t \in \mathcal{T}} \Bigg[ p^{\rm{im}}_t \Big(\lambda^{\rm{spot}}_t + Y^{\rm{im}}_t\Big) - p^{\rm{ex}}_t \Big(\lambda^{\rm{spot}}_t - Y^{\rm{ex}}_t\Big) \\
    & + (1-\beta^{\rm{DSO}}_t) Y^{\rm{im}}_t \Big( \sum_{i \in \mathcal{I}} p^+_{it} - p^{\rm{im}}_t \Big) + \alpha^{\rm{DSO}}_t p^{\rm{pen}}_t \Bigg], \label{eq:budget_balance}\
\end{align}    
where the left-hand side is the total payment of community members to the community manager during the day, whereas the right-hand side is identical to \eqref{eq:upper_obj} without the regularization, and corresponds to the total payment of the community operator for the electricity and network infrastructure usage.

Constraint \eqref{eq:individual_rationality} guarantees individual rationality for prosumers, meaning they are better off staying in the community rather than individually trading directly with the grid:
\begin{align}
    & \sum_{t \in \mathcal{T}} \Big( p^+_{it} - p^-_{it} \Big) x_{it} = C^{\, \rm{ext}}_i + w^+_i - w^-_i \quad \forall i \in \mathcal{I}, \label{eq:individual_rationality}
\end{align}  
where $w^+_{i},w^-_{i} \geq 0 \ \forall i \in \mathcal{I}$. 
The left-hand side calculates the total payment of prosumer $i$ as a community member under the dynamic pricing scheme. The right-hand side includes the pre-calculated value $C^{\, \rm{ext}}_i$,  corresponding to the payment of prosumer $i$ in a case under which the prosumer is not part of the community. If so, the prosumer pays the original grid tariffs without any discount and a proportional share of the penalty that would result from all prosumers optimizing their energy costs in an uncoordinated manner. This is all included in $C^{\, \rm{ext}}_i$. The right-hand side also contains two non-negative slack variables $w^+_i$ and $w^-_i$ to track how much better or worse off each prosumer is individually. These variables will later be used to formulate benefit distribution mechanisms among community members in Section \ref{sec:distribution}. In the case $w^-_i$ takes a positive value, it shows the underlying prosumer earns by being in the community, while a positive value of $w^+_i$ indicates the prosumer loses in the community. One may interpret $\sum_{i \in \mathcal{I}}w^-_i$ as the total benefit of forming the community (compared to a case where every prosumer acts individually). In contrast, $\sum_{i \in \mathcal{I}} w^+_i$ is the total cost if community formation is not beneficial. 
It is desirable that if a community member loses, no other member earns. This is enforced by
\begin{align}
    & \sum_{i \in \mathcal{I}} w^+_i \leq v^+, \ \sum_{i \in \mathcal{I}} w^-_i \leq v^-, \ \{ v^+, v^- \} \in \text{SOS1}, \label{eq:slack SOS1}
\end{align}  
where $v^+$ and $v^-$  are special order set of type 1 (SOS1) variables. By \eqref{eq:slack SOS1}, the total community benefit $\sum_{i \in \mathcal{I}}w^-_i$ and the total cost $\sum_{i \in \mathcal{I}} w^+_i$ cannot take positive values at the same time, i.e., at least one of these two terms should be zero.  

Constraint \eqref{eq:penalty} computes the amount of exceeded power $p^{\rm{pen}}_t$ imported by the community from the upstream grid, which is beyond the capacity limit $\bar{P}^{\rm{DSO}}_t$ imposed by the DSO:
\begin{align}
    0 \leq p^{\rm{pen}}_t &\geq p^{\rm{im}}_t - \bar{P}^{\rm{DSO}}_t \quad \forall t \in \mathcal{T}. \label{eq:penalty}\
    \end{align}  

The exported power is not penalized if it exceeds the cap to encourage local production injection to the grid. Constraint \eqref{eq:ul_reg} determines the maximum price $\bar{x}$ among prosumers and times, used for the regularization while also ensuring that no negative prices are set within the community:
\begin{align}
    0 \leq x_{it} &\leq \bar{x} \quad \forall i \in \mathcal{I}, \ t \in \mathcal{T}. \label{eq:ul_reg}\
\end{align}
    
The remaining upper-level constraints \eqref{eq:active power balance}-\eqref{eq:positive variables} are power flow constraints within the community. Constraints \eqref{eq:active power balance}-\eqref{eq:voltage limits} are the LinDistFlow equations, as originally proposed in \cite{Baran1989} and extensively applied in the literature, e.g., in \cite{Dvorkin2021,Mieth2018,Mieth2020,Berg2023}. It provides a lossless model for power flow in distribution grids. This set of constraints allows the community manager to ensure that the power flows at any point in time are feasible given the grid constraints of the energy community. Constraints \eqref{eq:active power balance} and \eqref{eq:reactive power balance} ensure that the total amount of active power and reactive power imported to or exported from the community at time $t$ is equal to the active power flow ($f^{\rm{p}}_{(n=0)t}$) or reactive power flow ($f^{\rm{q}}_{(n=0)t}$) at the interface, i.e.,  the reference node $n=0$:
%
\begin{align}
    p^{\rm{im}}_t - p^{\rm{ex}}_t &=  f^{\rm{p}}_{(n=0)t} \quad \forall t \in \mathcal{T} \label{eq:active power balance}\ \\
    q^{\rm{im}}_t - q^{\rm{ex}}_t &=  f^{\rm{q}}_{(n=0)t} \quad \forall t \in \mathcal{T}.  \label{eq:reactive power balance}\
\end{align}  

Constraint \eqref{eq:active power flow balance} calculates active power flow $f^{\rm{p}}_{nt}$ to every node $n$ except the reference node 0 at time $t$. Note that the set $i \in \mathcal{L}_n$ indicates prosumers $i$ located at node $n$, whereas the set $m \in \mathcal{D}_n$ refers to nodes $m$ downstream from node $n$. Constraint \eqref{eq:reactive power flow balance} is similar but for the reactive power flow $f^{\rm{q}}_{nt}$:
%
\begin{align}
    &f^{\rm{p}}_{nt} = \sum_{i \in \mathcal{L}_n}\Big( p^+_{it} - p^-_{it} \Big) + \sum_{m \in \mathcal{D}_n}  f^{\rm{p}}_{mt} && \forall t \in \mathcal{T}, \ n \in \mathcal{N}\setminus \{0\} \label{eq:active power flow balance}\\
    &f^{\rm{q}}_{nt} = \sum_{i \in \mathcal{L}_n}\Big( q^+_{it} - q^-_{it} \Big) + \sum_{m \in \mathcal{D}_n}  f^{\rm{q}}_{mt} && \forall t \in \mathcal{T}, \ n \in \mathcal{N}\setminus \{0\}. \label{eq:reactive power flow balance}\
\end{align}

Constraints \eqref{eq:base node}-\eqref{eq:voltage limits} bind nodal voltage magnitudes, where the squared voltage variable is replaced by an auxiliary variable $u_{nt}$ for modeling convenience. Constraint \eqref{eq:base node} sets the voltage magnitude of the reference node $n=0$ to 1 per unit. Constraint \eqref{eq:voltage} tracks the voltage drop throughout the nodes within the community, taking into account resistance and reactance parameters $R_n$ and $X_n$ for the line connecting node $n$ to the upstream one. Note that the set $m \in \mathcal{U}_{n}$ indicates all upstream nodes from node $n$. Finally, \eqref{eq:base node} constrains nodal voltage magnitudes to lie within $\underbar{$U$}$ and $\overline{U}$.

\begin{align}
    &u_{(n=0)t} = 1 \quad &&\forall t \in \mathcal{T} \label{eq:base node}\\
    &u_{nt} = \sum_{m \in \mathcal{U}_{n}}  u_{mt}  - 2 \Big( \frac{f^{\rm{p}}_{nt}}{S^{\mathrm{base}}}R_n + \frac{f^{\rm{q}}_{nt}}{S^{\mathrm{base}}}X_n\Big) \quad &&\forall t \in \mathcal{T}, \ n \in \mathcal{N}\setminus \{0\} \label{eq:voltage}\\
    & \underbar{$U$} \leq u_{nt} \leq \overline{U} \quad  &&\forall t \in \mathcal{T}, n \in \mathcal{N}\setminus \{0\}.\label{eq:voltage limits}\
\end{align}
    
Additionally, constraint \eqref{eq:line constraint} is added to enforce the apparent power flow capacity  $\overline{S}_n$ for the line that connects node $n$ to the upstream node in the radial distribution grid of the community:
    \begin{align}
    &(\frac{f^{\rm{p}}_{nt}}{S^{\mathrm{base}}})^2 + (\frac{f^{\rm{q}}_{nt}}{S^{\mathrm{base}}})^2 \leq \overline{S}_n^2 \quad &&\forall t \in \mathcal{T}, \ n \in \mathcal{N}\setminus \{0\}, \label{eq:line constraint}\ 
    \end{align}
which is a second-order cone constraint.

Constraints \eqref{eq:grid active power limit} and \eqref{eq:grid reactive power limit} enforce the physical capacity for the active and reactive power trade at the interface between the community and the upstream grid: 
    \begin{align}
    &p^{\rm{im}}_t,p^{\rm{ex}}_t  \leq \overline{P}^{\rm{grid}} \quad &&\forall t \in \mathcal{T} \label{eq:grid active power limit}\\
    &q^{\rm{im}}_t,q^{\rm{ex}}_t \leq \overline{Q}^{\rm{grid}} \quad &&\forall t \in \mathcal{T}. \label{eq:grid reactive power limit}\
        \end{align}

Finally, \eqref{eq:positive variables} declares the non-negativity conditions:
\begin{align}
&p^{\rm{im}}_t,p^{\rm{ex}}_t,q^{\rm{im}}_t,q^{\rm{ex}}_t \geq 0 \quad &&\forall t \in \mathcal{T}. \label{eq:positive variables}\
\end{align}
\end{subequations}

\subsubsection{Lower-level problem (optimal electricity dispatch)} \label{subsec:lower level}
We consider multiple prosumers in the energy community with flexible production, from rooftop photovoltaic panels and batteries, and/or inflexible consumption. The detailed formulation of the lower-level problems\footnote{We use, without loss of generality, the same generic mathematical formulation for different types of prosumers, and simply set the parameters and decision variables related to the operation of the photovoltaic panels and batteries to zero for prosumers who do not own such assets.} representing the individual optimal dispatch problem of each prosumer $i$ is given by \eqref{eqs:lower level}. The primal variable set for this problem is given by $\Phi_{i} = \{ p^{+}_{it},p^{-}_{it},q^{+}_{it},q^{-}_{it},p^{\rm{ch}}_{it},p^{\rm{dis}}_{it},e_{it},d^{\, \rm{shed}}_{it}\}$. Additionally, dual variables $\lambda_i^{(.)}$ and $\mu_i^{(.)}$ appear alongside their corresponding constraint as these are necessary when inserting the Karush-Kuhn-Tucker optimality conditions of the lower level into the upper-level optimization problem. Recall that $x_{it}$ is a parameter in the lower-level problem. The objective function of the lower-level problem, minimizing the cost of every prosumer $i$, is written as
\begin{subequations} \label{eqs:lower level}
    \begin{align}
        \underset{\Phi_i}{\operatorname{min}} \   \sum_{t \in \mathcal{T}} \Big[x_{it} \Big( p^+_{it} - p^-_{it} \Big)  + \alpha^{\rm{shed}} d^{\, \rm{shed}}_{it}\Big], \label{eq:lower obj}
    \end{align}
where $x_{it} \Big( p^+_{it} - p^-_{it} \Big)$ represents the payment that the prosumer  must make at time $t$ to the community manager for their energy use, whereas $\alpha^{\rm{shed}} d^{\, \rm{shed}}_{it}$ is the potential load shedding cost. Without loss of generality, we consider identical values of lost load $\alpha^{\rm{shed}}$ for all prosumers. 
Constraint \eqref{eq:ll power balance} is the power balance equation for each individual prosumer $i$:
    \begin{align}
        & p^+_{it} - p^-_{it} + PV_{it} - D_{it}  +d^{\, \rm{shed}}_{it}- p^{\rm{ch}}_{it} + p^{\rm{dis}}_{it} = 0 && \forall t \in \mathcal{T} &&: \lambda^{(1)}_{it},  \label{eq:ll power balance}
    \end{align}
where $PV_{it}$ and $D_{it}$ are deterministic solar production and demand forecasts, respectively. Each prosumer can choose to import or export power ($p^+_{it}$ and $p^-_{it}$,  respectively), charge or discharge their battery ($p^{\rm{ch}}_{it}$  and $p^{\rm{dis}}_{it}$), or shed power $d^{\, \rm{shed}}_{it}$.

Given battery charging and discharging efficiencies $\eta^{\rm{ch}}_i$ and $\eta^{\rm{dis}}_i$,  \eqref{eq:storage_1} ensures that the start and end level of the battery are the same, whereas  \eqref{eq:storage_all} is responsible for tracking the state of charge over the course of the day through the variable $e_{it}$:
    \begin{align}
        & e_{i(t=1)} = e_{i(t=24)} + \eta^{\rm{ch}}_i p^{\rm{ch}}_{i(t=1)} - \eta^{\rm{dis}}_i p^{\rm{dis}}_{i(t=1)} && &&:\lambda^{(2)}_{i} \label{eq:storage_1} \\
        & e_{it} = e_{i(t-1)} + \eta^{\rm{ch}}_i p^{\rm{ch}}_{it} - \eta^{\rm{dis}}_i p^{\rm{dis}}_{it} &&\forall t \in \mathcal{T}\setminus \{1\} &&:\lambda^{(3)}_{it}. \label{eq:storage_all}
    \end{align}

Assuming the reactive power demand is determined by a given relationship between the active and reactive power at the prosumer node, we enforce:
    \begin{align}
        & q^+_{it} = \sigma_i p^+_{it}  && \forall t \in \mathcal{T} &&: \lambda^{(4)}_{it} \label{eq:PQ positive} \\
        & q^-_{it} = \sigma_i p^-_{it}  && \forall t \in \mathcal{T} &&: \lambda^{(5)}_{it}, \label{eq:PQ negative}
    \end{align}
where $\sigma_i$ is given.    We declare non-negativity conditions as
    \begin{align}
        & p^+_{it}, p^-_{it}, q^+_{it}, q^-_{it}, p^{\rm{ch}}_{it}, p^{\rm{dis}}_{it}, e_{it}, d^{\, \rm{shed}}_{it} \geq 0  \ \ \forall t \in \mathcal{T} \notag\\ & \ \ \quad: \mu^{(1)}_{it}, \mu^{(2)}_{it}, \mu^{(3)}_{it}, \mu^{(4)}_{it}, \mu^{(5)}_{it}, \mu^{(6)}_{it},\mu^{(7)}_{it}, \mu^{(8)}_{it}. \label{eq:lower bound nobat}
    \end{align}
    
Finally, we enforce charging and discharging capacity $\bar{P}^{\rm{bat}}_i$ of the battery, energy storage capacity $\bar{E}_{i}$ of the battery, and upper limit of load curtailment:
    \begin{align}
        & p^{\rm{ch}}_{it}, p^{\rm{dis}}_{it} \leq \bar{P}^{\rm{bat}}_i && \forall t \in \mathcal{T} &&: \mu^{(9)}_{it}, \mu^{(10)}_{it} \label{eq:bat_variables_upper1}\\
        & e_{it} \leq \bar{E}_{i} && \forall t \in \mathcal{T} &&: \mu^{(11)}_{it} \label{eq:bat_variables_upper2}\\
        & d^{\, \rm{shed}}_{it} \leq D_{it} && \forall t \in \mathcal{T} &&: \mu^{(12)}_{it}. \label{eq:upper bound d} 
    \end{align}
\end{subequations}

As summary, the resulting bilevel problem is 
\begin{subequations} \label{eqs:summary}
\begin{align}
&\underset{\Omega}{\operatorname{min}} \quad  \eqref{eq:upper_obj} \\ 
& \text{s.t.} \quad \eqref{eq:budget_balance}-\eqref{eq:positive variables} \\ 
& p^{+}_{it}, p^{-}_{it}, q^{+}_{it}, q^{-}_{it}, d^{\, \rm{shed}}_{it} \in \underset{\Phi_i}{\text{argmin}} \ \Big\{ \eqref{eq:lower obj} \ \text{s.t.} \  \eqref{eq:ll power balance}-\eqref{eq:upper bound d} \Big\}.
\end{align}
\end{subequations}

Recall that the upper-level objective function \eqref{eq:upper_obj} is quadratic due to the price regularization term, whereas the upper-level constraints are linear, outside of the conic constraint \eqref{eq:line constraint}. Every lower-level problem is continuous and linear.

\subsection{Towards fair mechanisms for benefit allocation} \label{sec:distribution}

The total benefit $\sum_{i \in \mathcal{I}}w^-_i$ earned by forming a community to provide capacity limitation services (or in an extreme case, the total cost $\sum_{i \in \mathcal{I}}w^+_i$ if the community is unsuccessful) should be systematically distributed among community members. We introduce two distribution mechanisms, namely \textit{equal} and \textit{proportional} distributions.

Note that we do not consider any benefit distribution mechanism a posteriori. Instead, we propose a mechanism that considers fairness while setting dynamic prices. This is aligned with the proposal in \cite{zoe} that motivates fairness by design in shared-energy allocation problems. Both equal and proportional mechanisms proposed in this paper are based on adding an extra regularization term to the upper-level objective function \eqref{eq:upper_obj}, weighted by a positive weight $\gamma$.  Assigning a larger value for $\gamma$ further motivates the community manager to follow the underlying distribution but at the potential expense of a higher cost for the whole community. Both regularization terms preserve the convexity of \eqref{eq:upper_obj}. The regularization-based distribution mechanisms do \textit{not} guarantee that the targeted way of distributing the community benefit will be fully obtained. Rather, it \textit{motivates} the manager to assign prices accordingly. For example, in the case of the equal distribution mechanism, the total benefit may not be shared among members \textit{exactly} equally. However, it will be shared in a more equal way than a case without such a mechanism.

\subsubsection{Equal distribution} 
The equal distribution mechanism motivates an equal allocation of benefits among prosumers, irrespective of how much they consume or produce and how much flexibility they provide. This is achieved by penalizing \eqref{eq:upper_obj} by an extra regularization term, weighted by $\gamma^{\rm{eq}}$, such that the revised upper-level objective function can be written as
\begin{subequations}
    \begin{equation}
          \eqref{eq:upper_obj} +  \gamma^{\rm{eq}} \Bigg[\sum_{i\in \mathcal{I}} \Big( w^-_i -\widehat{w}^- \Big)^2 +  \sum_{i\in \mathcal{I}} \Big(  w^+_i  -\widehat{w}^+ \Big)^2 \Bigg],
     \end{equation}
where variables $\widehat{w}^-$ and $\widehat{w}^+$ are the mean values of $w_i^-$ and $w_i^+$ over $i\in \mathcal{I}$, respectively. Note that $\widehat{w}^-$ and $\widehat{w}^+$ should be added to the variable set $\Omega$. This regularizer penalizes any deviation from the mean benefit (if $\sum_{i \in \mathcal{I}}w^-_i$ takes a positive value) and the mean cost (if $\sum_{i \in \mathcal{I}}w^+_i$ takes a positive value). 

\subsubsection{Proportional distribution}
The proportional distribution mechanism motivates the community manager to allocate all prosumers a proportional amount of benefit (or cost) in line with their share of the total community baseline demand. Let the exogenous value $\Delta_{i}= \sum_t D_{it} - PV_{it}$ represent the total baseline residual load of prosumer $i$ over a day. We can define the demand share of prosumer $i$ as $\Delta_{i}/ \sum_i \Delta_{i}$. The proportional regularization term is then added to \eqref{eq:upper_obj}, weighted by  $\gamma^{\rm{pro}}$, such that the revised upper-level objective function reads as
%
%
\begin{align}
    \eqref{eq:upper_obj}  +   \gamma^{\rm{pro}}  \Bigg[ & \sum_{i\in \mathcal{I}} \Big( w^-_i  -    \frac{\Delta_i   }{\sum\limits_{i \in \mathcal{I}} \Delta_i} \sum_{i \in \mathcal{I}} w^-_i  \Big)^2 + \sum_{i\in \mathcal{I}} \Big( w^+_i  -    \frac{\Delta_i   }{\sum\limits_{i \in \mathcal{I}} \Delta_i}  \sum_{i \in \mathcal{I}} w^+_i  \Big)^2 \Bigg].
\end{align}
\end{subequations}
%

\vspace{-4mm}
\subsection{Resulting mixed-integer second-order cone program} 
\vspace{-1mm}
The following details the steps to reformulate the proposed bilevel program as a single-level mixed-integer second-order cone program (MISOCP). As the lower-level problems \eqref{eqs:lower level} are convex and satisfy Slater's conditions, their Karush-Kuhn-Tucker (KKT) conditions are both necessary and sufficient optimality conditions, and strong duality holds. Therefore, we replace the lower-level problems \eqref{eqs:lower level} by their equivalent KKT conditions in the upper level, and optimize over the set $\Theta = \big\{ \Omega, \Phi_i, \lambda_i^{(.)}, \mu_i^{(.)} \big\}$ which includes primal variables of the upper- and lower-level problems and dual variables of the lower-level problems. 
This results in a mathematical program with equilibrium constraints (MPEC). 

We then address two sources of non-convexity in this MPEC. First, we use the SOS1 reformulation as suggested by \cite{Siddiqui2013} and also applied in \cite{Askeland2023}, rather than the traditional Big-M method, to reformulate the complementarity conditions resulting from the lower-level inequality constraints, as the Big-M method can lead to incorrect solutions with poor M-value selections \cite{martin}. Second, as $x_{it}$, $p^+_{it}$, and $p^-_{it}$ are all decision variables of the resulting MPEC, \eqref{eq:budget_balance} and \eqref{eq:individual_rationality} contain bilinear terms  $x_{it}\left( p^+_{it} - p^-_{it} \right) $, which must be linearized. The application of the strong duality theorem to the lower-level problems allows for the replacement of such bilinear terms with an equivalent linear term from the lower-level dual problem, such that  
\begingroup
\allowdisplaybreaks
\begin{align} \label{eq:duality}
    \nonumber  \sum_{t \in \mathcal{T}} x_{it} \Big( p^+_{it} - p^-_{it} \Big) =& \sum_{t \in \mathcal{T}} \Bigg( \lambda^{(1)}_{it} \Big( PV_{it} - D_{it} \Big) - \bar{P}_{i}^{\rm{bat}} \Big(\mu^{(9)}_{it} + \mu^{(10)}_{it} \Big) \\
    & \quad \quad- \mu^{(12)}_{it}\bar{E}_{i} - \mu^{(12)}_{it} D_{it}  - \alpha^{\rm{shed}} d^{\, \rm{shed}}_{it} \Bigg) \ \ \ \forall{i} \in \mathcal{I}.
\end{align}
\endgroup

Finally, the quadratic regularization term $\rho \bar{x}^2$ in the upper-level objective function \eqref{eq:upper_obj} can be replaced by an auxiliary variable $z$. Then a second-order cone constraint $z \geq \rho \bar{x}^2$ can be added. 
The regularization terms for the benefit allocation can be treated similarly. By this, \eqref{eq:upper_obj} becomes linear, so the resulting problem will be a MISOCP, which can be solved by available commercial solvers for real-life applications.

\section{Numerical Study and Results} \label{sec:results}

\subsection{Setup for the case study} \label{subsec:case_study}

\begin{figure}[ht]
    \centering
    \includegraphics[width = 0.65\columnwidth]{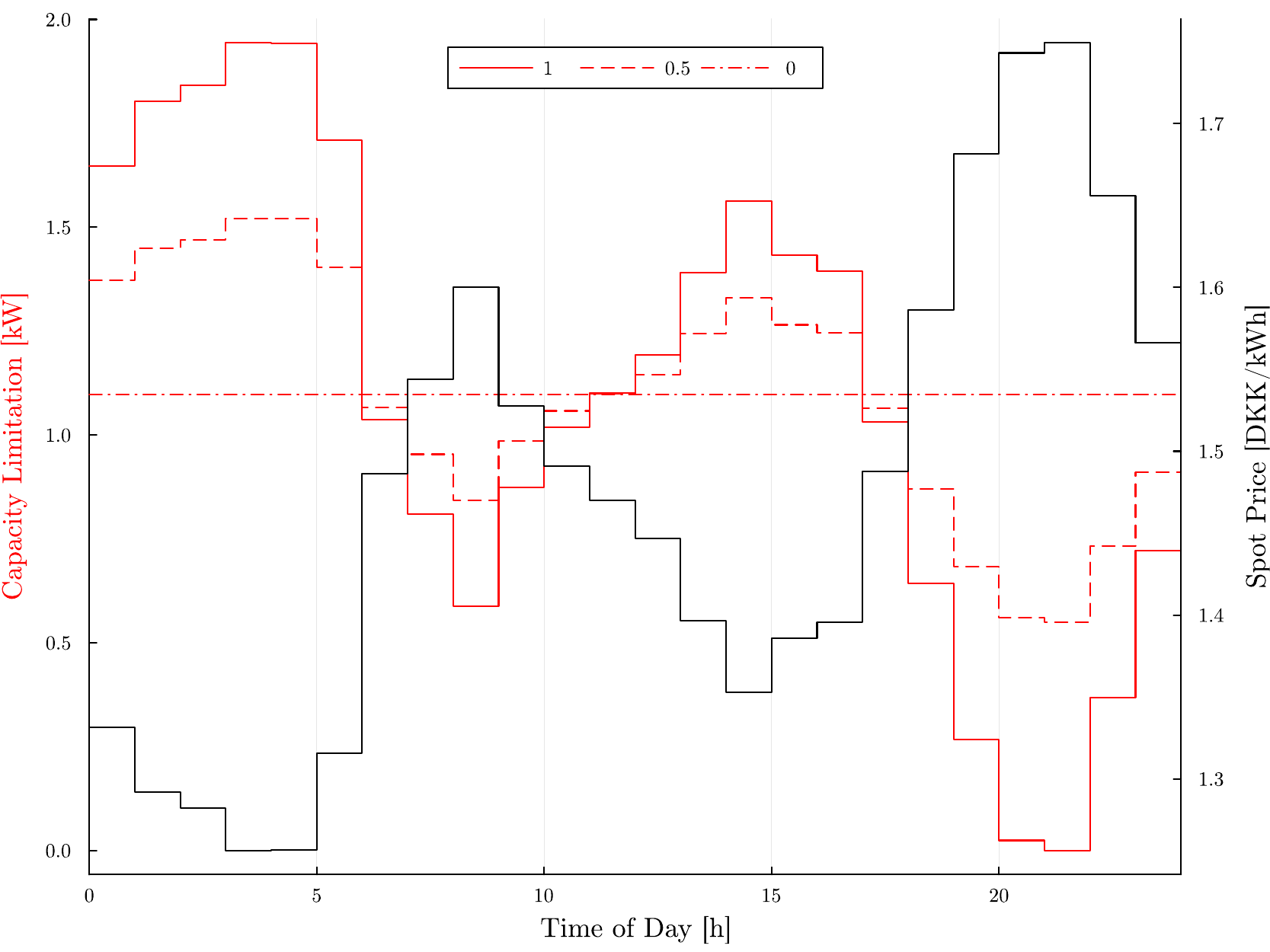}
    \caption{Three example curves for the capacity limitation $\bar{P}_t^{\rm{DSO}}$ over a day, with the capacity limitation variation factors of 0, 0.5, and 1 (red Y-axis). The DSO enforces either a constant or a time-variant capacity limit in an inversely proportional manner to the spot price (black Y-axis).}
    \label{fig:cap vs price}
\end{figure}

The proposed mechanism is applied to a 15-node radial distribution grid case study, as first introduced in \cite{Papavasiliou2018}, and additionally studied in \cite{Dvorkin2021, Mieth2018, Mieth2020}. We consider one prosumer per node, except the root node, each with a demand profile with an hourly resolution for a given day generated using the CREST model \cite{McKenna2020}. The apparent power base of the distribution system is chosen to scale these load profiles to similar per unit magnitude as used in \cite{Dvorkin2021}. Photovoltaic panels (PV) are distributed within the radial grid and have capacities ranging from 1 to 3 kW nominal capacity, each with a battery capable of storing 5 times the peak capacity of production. PV-battery pairs are placed at all nodes except nodes 3, 13, and 14. PV production profiles are generated at renewables.ninja using the methodology described in \cite{Pfenninger2016,Staffell2016}. Tariff values for the imported and exported community energy and contract penalty terms are based on existing Danish tariffs \cite{Radius2023,EnerginetTarifer2023} or chosen according to proposed Danish legislation \cite{KEF2022}. To discourage load shedding unless absolutely necessary due to grid congestion within the community, the value for load shedding is set to 25\% higher than the penalty for exceeding the DSO capacity limitation. The capacity limitation imposed by the DSO is inversely proportional to the day-ahead market price, as depicted in Figure \ref{fig:cap vs price}. The total available capacity under the capacity limitation is equal to the total residual load of the energy community. This reflects the increased demand in high-price hours, where a DSO might require a tighter capacity limitation to avoid congestion. To investigate the sensitivity of the community to this variation, a new parameter is introduced, the so-called \textit{capacity limitation variation factor}, varying from zero to one. This factor indicates to what extent $\bar{P}_t^{\rm{DSO}}$ will vary across hours. A value of 0 leads to a constant capacity limitation throughout the day, while a factor of 1 leads to a limitation with a minimum capacity of 0 kW and a higher maximum value. Figure \ref{fig:cap vs price} shows $\bar{P}_t^{\rm{DSO}}$ in kW (red Y-axis) when the capacity limitation variation factor is 0, 0.5, and 1. Figure \ref{fig:cap vs price} also illustrates the given spot market prices $\lambda^{\rm{spot}}_t$ in DKK/kW (black Y-axis). These prices are known to the community manager as the dynamic prices are set after the spot market has been cleared. Clearly, $\bar{P}_t^{\rm{DSO}}$ fluctuates more when the factor is 1, compared to the other cases.  Note that, for the sake of consistency, the total daily energy that can be imported by the community remains the same for all values of the capacity limitation variation factor. The corresponding code for this simulation is publicly available in \cite{Crowley2023}.

\begin{figure*}[t]
    \centering
    \includegraphics[width = \columnwidth]{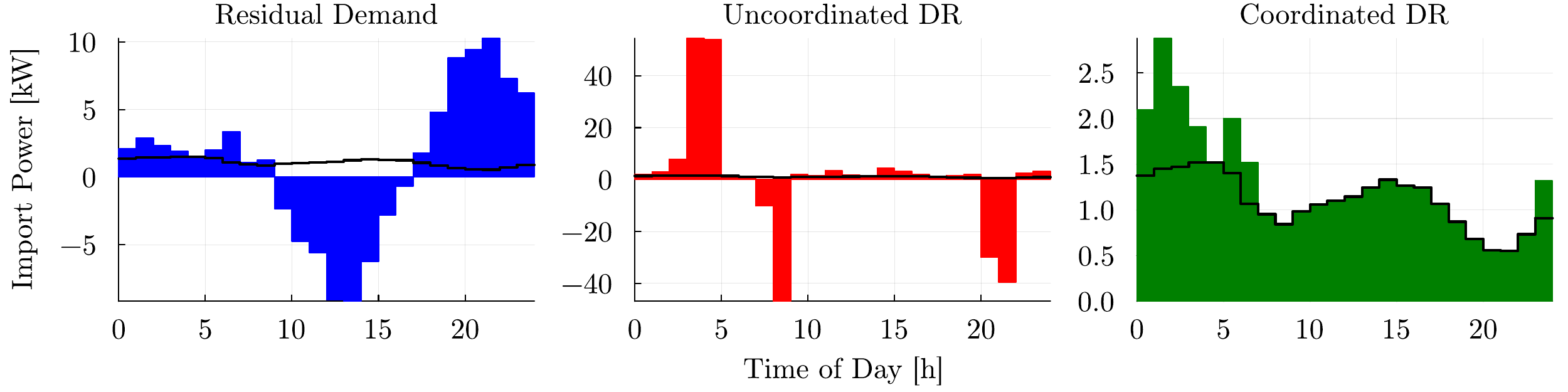}
    \caption{Power import (shaded areas) for different coordination mechanisms: with no community and no flexible prosumers (left); with prosumers shifting loads individually in an uncoordinated manner in response to spot prices (middle); and with a community manager coordinating the local flexibility of prosumers (right). The black curve represents the capacity limit $\bar{P}_t^{\rm{DSO}}$ imposed by the DSO. The abbreviation DR stands for demand response.}
    \label{fig:delivery}
\end{figure*}

\subsection{Delivery of capacity limitation services}

In this section, we investigate the benefits of the capacity limitation contract for the DSO, and in particular the success of the community in providing the contracted capacity limitation services. Firstly, the capacity limitation violation of the community members with and without a coordination mechanism are compared, as illustrated in Figure \ref{fig:delivery}. The left subplot corresponds to the base case, in which no community or flexibility is present in the distribution grid, and no temporal arbitrage is conducted in response to the spot prices. By this, the import power is simply the community's total residual demand (demand minus renewable production). In this case, the community load exceeds $\bar{P}_t^{\rm{DSO}}$ in some morning hours and all hours after 17:00. The middle subplot illustrates the so-called uncoordinated demand response (DR) case, in which there is no community and therefore flexible prosumers shift their loads in response to the spot market prices individually. As such, the prosumers respond without any awareness of upstream congestion or capacity limits imposed by the DSO. This naive load shifting causes a significant peak load at the lowest spot prices between 4:00 and 6:00, which would significantly exceed any DSO-imposed capacity limit as well as the apparent power limits of the distribution grid. Finally, the right subplot depicts the proposed coordinated DR mechanism, in which the community manager coordinates the response of the flexible prosumers through dynamic prices. A grid tariff discount factor of 0.5 and a capacity limitation variation factor of 1 is considered. This aims to investigate whether the community can adhere to even the most complex capacity limitation contract. With coordinated DR, the community still exceeds the capacity limitation in some hours, but by a significantly smaller margin than the other two cases. The reason the capacity limitation is exceeded, even though the capacity limitation curve has been designed to match the residual load of the community, is due to the charging and discharging inefficiencies of the batteries that are owned by the community members. Each time the battery is charged or discharged in an effort to meet the capacity limitation through temporal arbitrage, some energy is lost in the conversion process. This leads to an increase in the net load of the energy community, and hence the community exceeding the capacity limit. A future implementation of a capacity limitation contract for an energy community should account for the load increase brought upon by the temporal arbitrage, and increase the total available capacity throughout the day such that sufficient energy can be imported by the community to meet their demands and not penalize the community for trying to adhere to the capacity limit.

\begin{figure*}[t]
    \centering
    \includegraphics[width = \columnwidth]{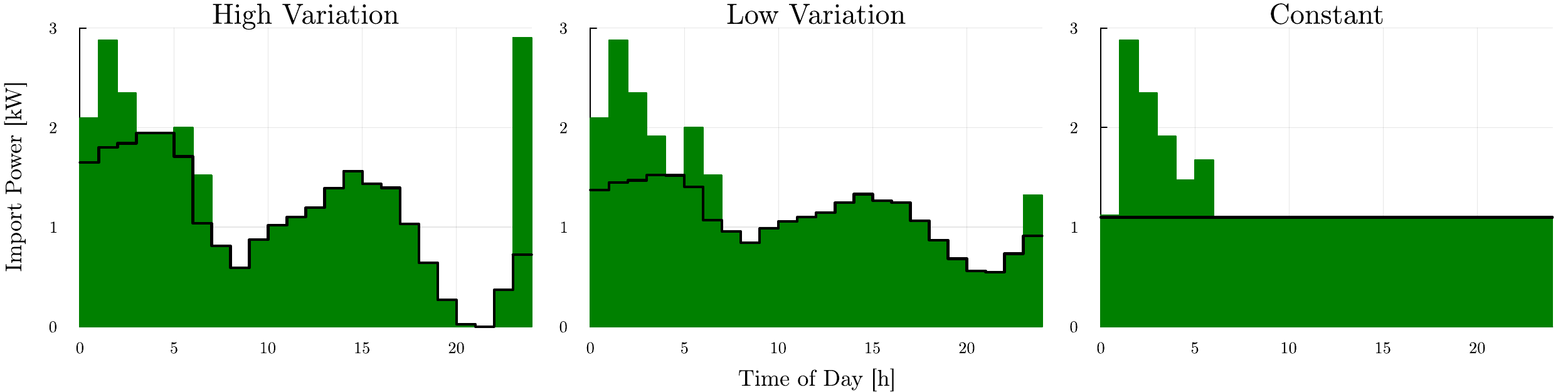}
            \vspace{-3mm}
    \caption{Power import of the community (shaded area) under capacity limitation variation factors of 1 (left), 0.5 (middle), and 0 (right). The black curve represents the capacity limit $\bar{P}_t^{\rm{DSO}}$ imposed by the DSO.}
    \label{fig:scenarios}
\end{figure*}

Furthermore, the impact of the variation factors on the energy community's delivery of capacity limitation services is investigated. The import profile of the community is illustrated in Figure \ref{fig:scenarios} for capacity limitation variation factors equal to 0 (constant), 0.5 (low variation), and 1 (high variation)\footnote{The left plot of Figure \ref{fig:scenarios} (high variation) is identical to the right plot of Figure \ref{fig:delivery}.}, while the discount tariff factor is kept constant at 0.5. In all three cases, the community violates the capacity limitation, with most violations occurring in the morning hours. As violating the capacity limitation has the same cost at all times, the deciding factor for when the community should exceed the capacity limit is the spot price and as the spot price is cheapest between 0:00 and 6:00, this is the cheapest time for the energy community to exceed the capacity limitation. The total violation is largest for the capacity variation factor of 1 because the batteries in the community are the most active trying to shift load, increasing the residual load the most. The largest single-hour violation occurs under this scenario as well. However, the single largest violation is in the evening. As the price in the last hour of the day is only slightly higher than in the morning, shifting this additional load to the morning would require too much additional battery charging and discharging and subsequent load increase, making a large violation in the evening the more cost-effective choice.

\subsection{Community cost}

\begin{figure}[t]
    \centering
    \includegraphics[width = 0.65\columnwidth]{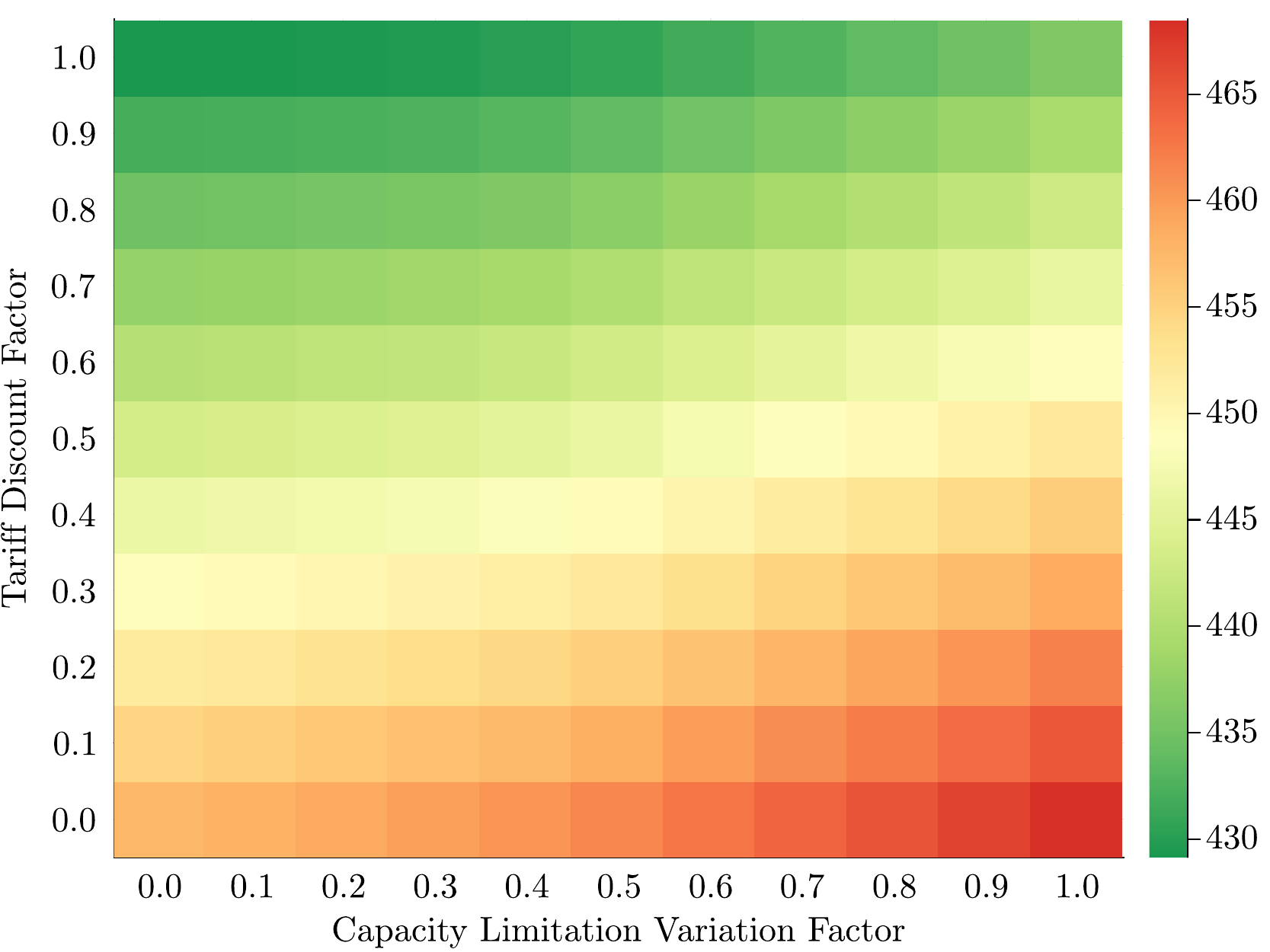}
            \caption{Community cost (in DKK) over a day as a function of the tariff discount factor $\beta^{\rm{DSO}}_t$ and the capacity limitation variation factor. Note that all cases are beneficial compared to the Uncoordinated DR case shown in Figure \ref{fig:delivery}. }
    \label{fig:heatmap}
\end{figure}

In this section, we investigate the benefits of the capacity limitation contract for the community as a whole. In particular, the impact of the discounted DSO tariff ($1-\beta^{\rm{DSO}}_t$) and the capacity limitation variation factor on the total community cost is illustrated in Figure \ref{fig:heatmap}. Recall that $\beta^{\rm{DSO}}_t$ ranges between zero (full tariff) to one (no tariff) and applies for a whole day. Similarly, the capacity limitation variation factor takes values between zero (no variation) and one (high variation). As observed along the y-axis of Figure \ref{fig:heatmap}, the overall cost decreases when the community pays a lower percentage of the DSO grid tariff for their internal power flows. This relationship is linear, and for every 10\% that the tariff discount factor increases, the savings for the entire energy community increase by 3.2 DKK (0.43€) on average. In contrast, as observed along the x-axis, the community cost also increases when the capacity limitation variation factor increases. This relationship can be explained by the increased use of the batteries required to shift loads between hours, leading to increased charging and discharging losses and, as a result, higher total energy imported and penalties. This phenomenon is slightly offset by the ability of the community to import more energy during morning hours with low spot prices when the capacity limitation is larger. However, this is insufficient to make a higher capacity limitation variation factor beneficial for the community. We observe that the relationship between the community cost and capacity limitation variation factor is exponential. Between capacity limitation variation factors of 0.0 and 0.3, the change in cost is minimal, but the growth rate of the community cost increases for higher values of the capacity limitation variation factor.


\subsection{Fairness, allocation of community benefits, and dynamic pricing}

In this section, we investigate the benefits of the capacity limitation contract for individual prosumers and fairness within the community. All results in this section consider a specific capacity limitation contract with a variation factor of 1 and a tariff discount factor of 0.5. Firstly, the distribution of the community benefit among the prosumers under two distribution mechanisms is investigated. The community benefit refers to $\sum_{i \in \mathcal{I}}w^-_i$,  i.e., the collective savings of the community with respect to the uncoordinated DR case under which there is no community and all prosumers act individually. The import profile for this benchmark case is provided in the middle subplot of Figure \ref{fig:delivery}. It is assumed that in the uncoordinated DR case, prosumers pay for the costs incurred by the significant peak demand in hours 4 and 5 at the same penalty rate as under the coordinated DR case (75 DKK/kW). A small weight $\gamma=10^{-6}$ is assigned to the regularization term for each distribution mechanism. This value is chosen as the community's import and export power profiles are identical to the case with no distribution mechanism under this weight, even though the prices for individual consumers have changed.

\begin{figure}[t]
    \centering
    \includegraphics[width = 0.65\columnwidth]{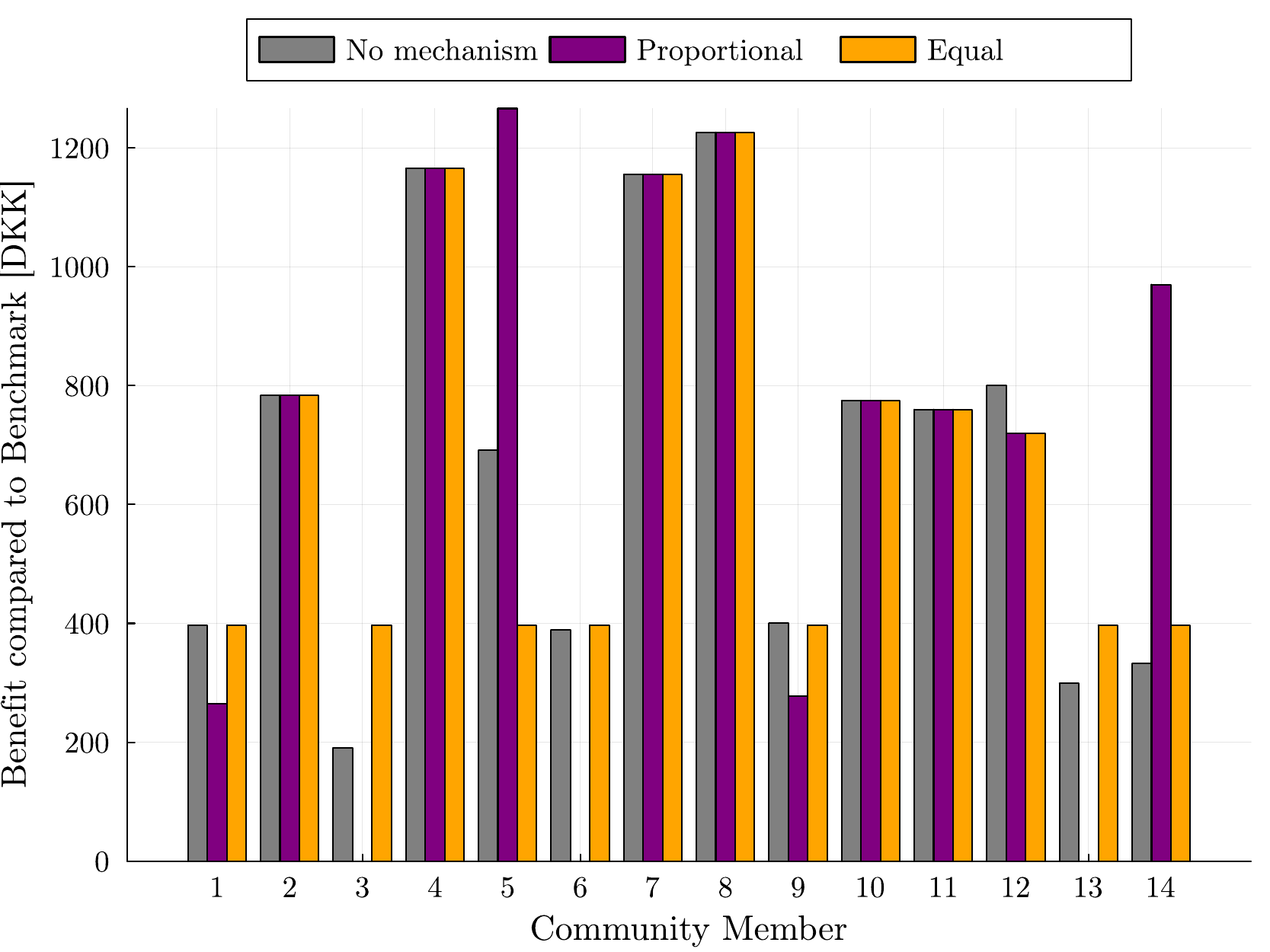}
        \caption{Invidual benefit (in DKK) of each community member under three distribution mechanisms (none, equal, and proportional).}
    \label{fig:distribution}
\end{figure}

Figure \ref{fig:distribution} shows how the community benefit is shared among the 14 community members without a distribution mechanism and when the community manager aims for a proportional or equal distribution. In the case of the proportional distribution, a similar variation of benefit distribution is observed among the prosumers as in the case with no distribution. However, as intended, the mechanism tends to distribute the highest shares of the community benefits to the members with the highest residual demand (prosumers 5 and 14). On the other hand, the equal distribution mechanism ensures that prosumers 1, 3, 5, 6, 9, 13, and 14 all earn equal amounts of benefit, but this is not achieved for other members of the community. In addition, the equal distribution mechanism reduces the maximum benefit redistributed among the prosumers compared to the proportional distribution mechanism. For some members of the community (prosumers 2, 4, 7, 8, 10 and 11) there is no observable difference in distributed benefit between various distribution mechanisms. It is worth noting that all members who observe no change in benefit distribution have a negative residual load. The benefit distribution mechanisms do not seem capable of addressing this in their current formulation. As previously stated, the distribution mechanisms are implemented through a regularization term with a small weight in the objective function. Therefore equality or proportionality of benefit cannot be guaranteed. A higher regularization weight would prioritize the fairness term in the objective function, over the community cost term. However, this would come at the expense of optimality in delivering capacity limitation services, which remains the overarching goal of the proposed mechanism in this paper.

\begin{figure}[t]
    \centering
    \includegraphics[width=\linewidth]{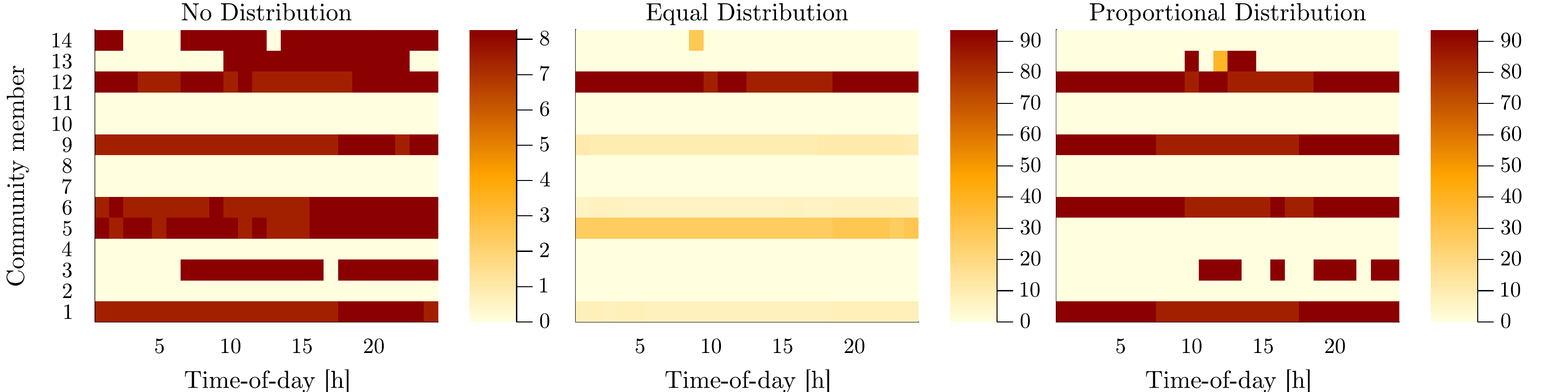}
    \caption{Heatmaps of individual prices [DKK/kWh] for all prosumers throughout the day under three benefit distribution mechanisms (none, equal, and proportional)}
    \label{fig:heatmap prices}
\end{figure}

Besides, the impact of the distribution mechanisms on the dynamic prices $x_{it}$ for each community member is investigated. Figure \ref{fig:heatmap prices} shows three heatmaps of the dynamic prices for each community member over a day with and without a distribution mechanism. We observe that the scale of the prices with no distribution mechanism is significantly lower compared to the scale of the prices with equal or proportional distribution mechanisms. This shows that without an attempt at fair distribution, the community manager sets comparatively low prices for all members of the community that achieve the desired response to satisfy the capacity limitation contract. Both distribution mechanisms drastically increase prices as a whole but also change individual price patterns. The maximum price under both distribution mechanisms is limited by the load-shedding cost for all prosumers, equal to 93.75 DKK/kWh. Tuning the weights of the benefit regularization terms and the maximum price regularizer may prove beneficial in finding a middle ground and not raising the prices directly to the highest possible amount. 

A large portion of the community members observe a price of 0 DKK/kWh throughout the day or at different times. For example, prosumers 2, 4, 7, 8, 10, and 11 have a price of 0 DKK/kWh in all hours under all distribution mechanisms. These prosumers have a negative residual load over the day and do not have to pay for energy consumption. Therefore the remaining subset of prosumers must cover the cost of the imported electricity, with the regularization mechanisms impacting who pays the least and who pays the most. Furthermore, the price curves throughout all three distribution cases are quite volatile, regularly setting the prices at the maximum or minimum; in fact, all price curves have a maximum of three prices throughout the day (minimum, maximum, intermediate). Such a price curve differs significantly from the smoother cost curves obtained in Lagrange-multiplier-based day-ahead markets. This difference can be attributed to the limiting modeling choices for the prosumers' price response, modeled as simple linear problems. In such a linear problem, optimal responses switch between vertices of their feasible region at a certain price threshold and are indifferent to a wide range of prices in between. A more nuanced prosumer model might prove beneficial in diversifying the price signatures and allowing for smoother cost curves and optimal responses.

\subsection{Computational scalability} \label{subsec:results_comp}

In this section, the scalability and computational tractability of the proposed mechanism are investigated. The proposed bilevel model is solved on three additional expanded case studies. The same case study described in Section \ref{subsec:case_study} is used, but the grid is extended by adding the same radial grid structure at node 14 to increase the number of nodes and prosumers. This is done multiple times to create case studies scaled up by a factor of 2, 4, and 8, respectively. All results in this section consider a specific capacity limitation contract with a variation factor of 1 and a tariff discount factor of 0.5. Similarly to the approach in \cite{Berg2023}, for the sake of computational tractability, the bilevel model is solved in two steps: \textit{i}) the model is solved without including the second-order conic constraint \eqref{eq:line constraint}; and \textit{ii}) an ex-post check of power flows and line capacities is conducted to ensure the feasibility of the solution. We observe that all simulations satisfy constraint \eqref{eq:line constraint}. Table \ref{tab:computational} summarizes the number of variables, SOS1 constraints, and the computational time of each simulation. A 2020 16-GB Apple M1 MacBook Pro is used for simulations. With each doubling of nodes in the case study, all three factors roughly double in scale, suggesting that scaling this model should not prove problematic with this 2-step approach.

\begin{table}[t]
    \centering
    \begin{tabular}{|ccc|c|}
        \hline
        \# of Members & \# of Variables & \# of SOS1 constraints & \multicolumn{1}{|c|}{Computational time (s)}\\ \hline
        14 & 21,679 & 4,033 & 0.31 \\
        28 & 43,211 & 8,065 & 0.72 \\
        56 & 86,275 & 16,129 & 2.67 \\ 
        112 & 172,403 & 32,257 & 9.40 \\ \hline
    \end{tabular}
    \caption{Analysis of the computational scalability of the proposed bilevel model.}
    \label{tab:computational}
\end{table}

\section{Conclusion} \label{sec:conclusion}

This paper shows how an energy community can meet the requirements of a capacity limitation service contract based on dynamic pricing, for which we develop a mathematical framework based on bilevel programming. This framework enables the provision of capacity limitation services by the community to the DSO while ensuring desirable market properties, including budget balance for the community manager and individual rationality for every community member. A regularization-based mechanism encourages a fair allocation of collected benefits among community members. 

A numerical study is conducted to show the success of the proposed community framework compared to a case in which there is no community and prosumers respond individually to spot prices in an uncoordinated manner. A few lessons can be learned from the numerical study that could be useful for DSOs and future community managers. Firstly, the community benefits the most from an invariant capacity limitation as this requires the least charging and discharging of the battery, thereby increasing the net load by the smallest amount. Secondly, the increased net load due to temporal arbitrage using batteries must be considered when setting a capacity limitation for the energy community. Not doing so would unfairly penalize the community for trying to adhere to the DSO capacity limit. Lastly, by dropping the second-order conic constraint for the apparent power flow limits, which is unlikely to be constraining the optimization in the first place, the bilevel pricing model can be easily scaled up to handle large-scale energy communities, potentially providing larger synergies within the community and the DSO.

Future work in this framework should extend the proposed mechanism to account for stochasticity and non-stationarity in the community members' demand profiles and local production. This work looks to provide an ideal deterministic benchmark with an optimistic view of lower-level problems, showing the full potential of a community in delivering capacity limitation services. Therefore, the mechanism relies on accurate forecasting and truthful reporting from community members to price the electricity correctly, but a stochastic extension would lead to a more practical price-setting process. A more realistic implementation would require a more diverse asset pool and consideration for more nuanced consumer behavior. In the numerical study, the prosumers exposed to high prices still prefer staying in the community over the load-shedding alternative. However, rationality among community members cannot be guaranteed in a real-life setting, and more risk-averse community members may not be pleased with such volatility. Including a more refined model of consumer behavior in lower-level problems and learning the demand response behaviors of various community members over time could be beneficial. One could explore the applications of multi-armed bandit or similar learning approaches to achieve this. It would be interesting to explore how to derive dynamic prices in this setup as well as the regret faced by the community manager as they learn about community members' behavior and whether desirable market properties can still be guaranteed.

\section*{Acknowledgement} This research was supported by the Danish Energy Technology Development and Demonstration Programme (EUDP) through the Flexible Citizen Energy Communities (FLEX-CEC) project (64021-1090).

\appendix

\end{document}